\documentclass[11pt]{article}
\newcommand		{\comment}[1]		{}

\usepackage{amsmath,amsthm} 			
\usepackage{mathtools}

\usepackage[left=2.54cm,right=2.54cm,top=2.54cm,bottom=2.54cm]{geometry}
\usepackage[nouppercase]{scrlayer-scrpage}
\usepackage{xspace}
\usepackage{enumitem}

\usepackage{avant} 				
\usepackage[sc,osf]{mathpazo} 			
	\AtBeginDocument{
		\DeclareSymbolFont{AMSb}{U}{msb}{m}{n}
		\DeclareSymbolFontAlphabet{\mathbb}{AMSb}
	}
\usepackage{extpfeil} 
\usepackage{mathabx}  
\usepackage{mathrsfs}	
\usepackage{fancyvrb}	

\usepackage{rotating,pdflscape}
\usepackage{afterpage}
\usepackage{graphicx}
\usepackage{float}
\usepackage{subcaption}
\usepackage{caption}
\newcommand{\mockalph}[1]{\!}

\usepackage[utf8]{inputenc}
\usepackage[T1]{fontenc}

\usepackage{tocloft}
\makeatletter
\renewcommand{\l@figure}{\@dottedtocline{1}{1em}{3.5em}}
\renewcommand{\l@table}{\@dottedtocline{2}{1em}{3.5em}}
\makeatother
\newcommand*{\noaddvspace}{\renewcommand*{\addvspace}[1]{}}
\addtocontents{lof}{\protect\noaddvspace}

\usepackage[hyphens]{url}		%
\usepackage{hyperref}	
\usepackage{cleveref} 	
\newcommand		{\myred}		{BrickRed}
\hypersetup{
			colorlinks		= true, 	
			urlcolor		= \myred, 
			linkcolor		= \myred, 	
			citecolor		= PineGreen 	
}
\newcommand		{\hyref}[1]		{\hyperref[#1]{\ref*{#1}}}

\newif\ifdebug                                                      %
\debugfalse
\ifdebug\usepackage{lineno}\linenumbers\else\fi

\usepackage{etoolbox}
\PassOptionsToPackage{dvipsnames,svgnames,x11names,table}{xcolor}
\newcommand		{\defd}[1]	{\textcolor{RoyalBlue}{\textbf{\textit{#1}}}}
\newcommand		{\defm}[1]	{\textcolor{RoyalBlue}{#1}}

\usepackage[all]{xypic}\xyoption{rotate}
	\newdir{ >}{{}*!/-7pt/\dir{>}}
	\newdir{i}{{}*!/-7pt/\dir{(}	}
\usepackage{tikz}
\usetikzlibrary{matrix,fit,backgrounds,calc,arrows} 
\tikzstyle{image}=[rectangle,fill=Red!20,inner sep=-2pt]
\tikzstyle{nonzero}=[rectangle,fill=Navy!20,inner sep=0pt]
\tikzstyle{nonzerosm}=[rectangle,fill=Navy!20,inner sep=-2pt]

\usepackage{thmtools}
\makeatletter
\let\c@figure\c@table
\let\c@equation\c@table
\makeatother

\numberwithin{table}{section}
\numberwithin{figure}{section}
\newtheorem{abcthm}[table]{Theorem}

\newtheorem{theorem}[table]{Theorem}

\newtheorem{proposition}[table]{Proposition}

\newtheorem{corollary}[table]{Corollary}

\newtheorem{lemma}[table]{Lemma}

\newtheorem{claim}[table]{Claim}

\theoremstyle{definition}

\newtheorem{definition}[table]{Definition}

\newtheorem{notation}[table]{Notation}

\newtheorem{observation}[table]{Observation}

\newtheorem{conjecture}[table]{Conjecture}

\theoremstyle{remark}
\newtheorem{fact}[table]{Fact}

\newtheorem{example}[table]{Example}

\newtheorem{exercise}[table]{Exercise}

\newtheorem{problem}[table]{Problem}
\newtheorem{histrmks}[table]{Historical remarks}
\newtheorem{remark}[table]{Remark}
\newtheorem{remarks}[table]{Remarks}

\theoremstyle{plain}
\newtheorem*{thm*}{Theorem}
\newtheorem*{theorem*}{Theorem}
\newtheorem*{prop*}{Proposition}
\newtheorem*{proposition*}{Proposition}
\newtheorem*{lemma*}{Lemma}
\newtheorem*{corollary*}{Corollary}
\newtheorem*{cor*}{Corollary}

\theoremstyle{definition}
\newtheorem*{definition*}{Definition}
\newtheorem*{defn*}{Definition}
\newtheorem*{QQ*}{Question}
\newtheorem*{obs*}{Observation}
\newtheorem*{notation*}{Notation}
\newtheorem*{discussion*}{Discussion}

\theoremstyle{remark}
\newtheorem*{rmk*}{Remark}
\newtheorem*{remark*}{Remark}
\newtheorem*{examples*}{Examples}
\newtheorem*{example*}{Example}
\newtheorem*{EG*}{Example}
\newtheorem*{EGs*}{Examples}
\newtheorem*{fact*}{Fact}
\newtheorem*{prob*}{Problem}

\newcommand{\bthm}{\begin{theorem}}
\newcommand{\ethm}{\end{theorem}}
\newcommand{\bprop}{\begin{proposition}}
\newcommand{\eprop}{\end{proposition}}
\newcommand{\bcor}{\begin{corollary}}
\newcommand{\ecor}{\end{corollary}}
\newcommand{\bconj}{\begin{conjecture}}
\newcommand{\econj}{\end{conjecture}}
\newcommand{\blem}{\begin{lemma}}
\newcommand{\elem}{\end{lemma}}
\newcommand{\bclm}{\begin{claim}}
\newcommand{\eclm}{\end{claim}}
\newcommand{\bpf}{\begin{proof}}
\newcommand{\epf}{\end{proof}}
\newcommand{\bdetails}{\begin{details}}
\newcommand{\edetails}{\end{details}}
\newcommand{\bdefi}{\begin{definition}}
\newcommand{\edefi}{\end{definition}}
\newcommand{\bdefn}{\begin{definition}}
\newcommand{\edefn}{\end{definition}}
\newcommand{\bex}{\begin{example}}
\newcommand{\eex}{\end{example}}
\newcommand{\bprob}{\begin{problem}}
\newcommand{\eprob}{\end{problem}}
\newcommand{\bob}{\begin{observation}}
\newcommand{\eob}{\end{observation}}
\newcommand{\bexer}{\begin{exercise}}
\newcommand{\eexer}{\end{exercise}}
\newcommand{\bexers}{\begin{exercises}}
\newcommand{\eexers}{\end{exercises}}
\newcommand{\brmk}{\begin{remark}}
\newcommand{\ermk}{\end{remark}}
\newcommand{\bhist}{\begin{histrmks}}
\newcommand{\ehist}{\end{histrmks}}
\newcommand{\brmks}{\begin{remarks}}
\newcommand{\ermks}{\end{remarks}}
\newcommand{\bverb}{\begin{verbatim}}
\newcommand{\everb}{\end{verbatim}}
\newcommand{\bntn}{\begin{notation}}
\newcommand{\entn}{\end{notation}}
\newcommand{\bfct}{\begin{fact}}
\newcommand{\efct}{\end{fact}}
\newcommand{\bfcts}{\begin{facts}}
\newcommand{\efcts}{\end{facts}}
\newcommand{\benum}{\begin{enumerate}}
\newcommand{\eenum}{\end{enumerate}}
\newcommand{\bitem}{\begin{itemize}}
\newcommand{\eitem}{\end{itemize}}

\tracingpatches
\makeatletter
\patchcmd{\@setref}{\bfseries ??}{\bfseries\color{red} FIX ME!}{}{}
\patchcmd{\@setcite}{\bfseries ?}{\bfseries\color{red} FIX ME!}{}{}
\patchcmd{\@setcref}         {??}{\color{red} FIX ME!}{}{}
\patchcmd{\@setcref}         {??}{\color{red} FIX ME!}{}{}
\patchcmd{\@setcrefrange}    {??}{\color{red} FIX ME!}{}{}
\patchcmd{\@setcrefrange}    {??}{\color{red} FIX ME!}{}{}
\patchcmd{\@setcrefrange}    {??}{\color{red} FIX ME!}{}{}
\patchcmd{\@setcrefrange}    {??}{\color{red} FIX ME!}{}{}
\patchcmd{\@setcrefrange}    {??}{\color{red} FIX ME!}{}{}
\patchcmd{\@setcrefrange}    {??}{\color{red} FIX ME!}{}{}
\patchcmd{\@setnamecref}     {??}{\color{red} FIX ME!}{}{}
\patchcmd{\@setnamecref}     {??}{\color{red} FIX ME!}{}{}
\patchcmd{\@setcpageref}     {??}{\color{red} FIX ME!}{}{}
\patchcmd{\@setcpageref}     {??}{\color{red} FIX ME!}{}{}
\patchcmd{\@setcpagerefrange}{??}{\color{red} FIX ME!}{}{}
\patchcmd{\@setcpagerefrange}{??}{\color{red} FIX ME!}{}{}
\patchcmd{\@setcpagerefrange}{??}{\color{red} FIX ME!}{}{}
\patchcmd{\@setcpagerefrange}{??}{\color{red} FIX ME!}{}{}
\patchcmd{\@setcpagerefrange}{??}{\color{red} FIX ME!}{}{}
\patchcmd{\@cref}            {??}{\color{red} FIX ME!}{}{}
\def\blx@citation@entry#1#2{%
  \blx@bibreq{#1}%
  \ifinlist{#1}{\blx@cites}
    {}
    {\listgadd{\blx@cites}{#1}%
     \blx@auxwrite\@mainaux{}{\string\abx@aux@cite{#1}}}%
  \ifinlistcs{#1}{blx@segm@\the\c@refsection @\the\c@refsegment}
    {}
    {\listcsgadd{blx@segm@\the\c@refsection @\the\c@refsegment}{#1}}%
  \blx@ifdata{#1}%
    {}%
    {\ifcsdef{blx@miss@\the\c@refsection}%
       {\ifinlistcs{#1}{blx@miss@\the\c@refsection}%
          {{\bfseries\color{red} cite:} }%
          {\blx@logreq@active{#2{#1}}}}%
       {\blx@logreq@active{#2{#1}}}}}
\def\blx@citeadd#1{%
  \ifcsdef{blx@keyalias@\the\c@refsection @#1}
    {\edef\blx@realkey{\csuse{blx@keyalias@\the\c@refsection @#1}}}
    {\def\blx@realkey{#1}}%
  \expandafter\blx@citation\expandafter{\blx@realkey}\blx@msg@cundefon
  \expandafter\blx@ifdata\expandafter{\blx@realkey}
    {\advance\blx@tempcnta\@ne
     \listeadd\blx@tempa{\blx@realkey}}
    {\ifnum\blx@tempcntb>\z@\multicitedelim\fi
     \expandafter\abx@missing\expandafter{\blx@realkey}%
     \advance\blx@tempcntb\@ne}}
\makeatother

\newcommand{\presectionskip}{-1.5\baselineskip}
\newcommand{\postsectionskip}{0.3\baselineskip}
\makeatletter
\usepackage{titlesec}
\renewcommand{\section}{\@startsection
  {chapter}{0}{0mm}
  {\presectionskip}
  {\postsectionskip}
  {\sffamily\huge}}
\renewcommand{\section}{\@startsection
  {section}{1}{0mm}
  {\presectionskip}
  {\postsectionskip}
  {\sffamily\LARGE}}
\renewcommand{\subsection}{\@startsection
  {subsection}{2}{0mm}
  {\presectionskip}
  {\postsectionskip}
  {\sffamily\Large}}
\renewcommand{\subsubsection}{\@startsection
  {subsubsection}{3}{0mm}
  {\presectionskip}
  {\postsectionskip}
  {\sffamily\normalsize}}
\renewcommand{\@seccntformat}[1]{\csname the#1\endcsname.\quad}
\newcommand\HUGE{\@setfontsize\Huge{30}{47}} 
\makeatother
\usepackage{titling}
  \titleformat{\chapter}[display]
  {\sffamily\Large}
  {Chapter {\HUGE\normalfont\thechapter}}    
  {1em}
  {\huge}
  \titleformat{\section}
  {\sffamily\Large}
  {\normalfont{\@setfontsize\Large{18}{47}\thesection.}}    
  {1em}
  {}
  \titleformat{\subsection}
  {\sffamily\large}
  {\normalfont{\Large\thesubsection}.}    
  {1em}
  {}
  \titleformat{\subsubsection}
  {\sffamily\normalsize}
  {\normalfont{\large\thesubsubsection}.}    
  {1em}
  {}

\renewcommand		{\SS}				{\textsection}

\newcommand		{\quation}[1]			{\begin{equation} #1 \end{equation}}

\newcommand		{\eqn}[1]			{\begin{align*} #1 \end{align*}}

\def			\SPSB#1#2			{\rlap{\textsuperscript{#1}}\textsubscript{#2}}
\makeatletter
\def			\smallunderbrace#1		{\mathop{\vtop{\m@th\ialign{##\crcr
							   $\hfil\displaystyle{#1}\hfil$\crcr
							   \noalign{\kern3\p@\nointerlineskip}%
							   \tiny\upbracefill\crcr\noalign{\kern3\p@}}}}\limits}
\makeatother

\makeatletter
\newcommand{\subalign}[1]{%
  \vcenter{%
    \Let@ \restore@math@cr \default@tag
    \baselineskip\fontdimen10 \scriptfont\tw@
    \advance\baselineskip\fontdimen12 \scriptfont\tw@
    \lineskip\thr@@\fontdimen8 \scriptfont\thr@@
    \lineskiplimit\lineskip
    \ialign{\hfil$\m@th\scriptstyle##$&$\m@th\scriptstyle{}##$\crcr
      #1\crcr
    }%
  }
}
\makeatother

\makeatletter
\newcommand		{\oset}[3][0ex]			{%
								\raisebox{.175ex}{$%
								  \mathrel{\mathop{#3}\limits^{
								    \vbox to#1{\kern-2\ex@
								    \hbox{$\scriptstyle#2$}\vss}}}
								    $}%
							    }
\makeatother

\usepackage{faktor,xfrac}

\newcommand		{\dsp}		{\displaystyle}
\newcommand		{\nd}		{\noindent}

\newcommand		{\bs}		{\bigskip}

\newcommand		{\mn}		{\mspace{-2mu}}
\newcommand		{\mnn}		{\mspace{-1mu}}

\newcommand		{\ol}			{\overline}
\newcommand		{\os}			{\overset}
\newcommand		{\us}			{\underset}

\newcommand		{\wt}			{\widetilde}

\newcommand		{\mr}			{\mathrm}
\newcommand		{\bb}			{\mathbb}

\newcommand		{\ms}			{\mathscr}
\newcommand		{\sans}			{\mathsf}

\newcommand		{\g}		{\gamma}

\renewcommand		{\epsilon}	{\varepsilon}
\newcommand		{\e}		{\epsilon}

\newcommand		{\h}		{\eta}

\renewcommand		{\l}		{\lambda}

\newcommand		{\s}		{\sigma}

\newcommand		{\W}		{\Omega}

\newcommand		{\D}		{\Delta}

\DeclareSymbolFont{cmletters}{OT1}{cmr}{m}{n}
\DeclareMathSymbol{\Ups}{\mathalpha}{cmletters}{"7}
\renewcommand		{\Upsilon}	{\Ups}

\newcommand		{\ceq}		{\coloneqq}

\renewcommand		{\cup}		{\mspace{-1mu}\smile\mspace{-1mu}}
\makeatletter
\DeclarePairedDelimiterX
			{\pmodx}[1]	{(}{)}{{\operator@font mod}\mkern6mu#1}
					\renewcommand{\pmod}{%
					  \allowbreak
					  \if@display\mkern18mu\else\mkern8mu\fi
						  \pmodx
					}
\makeatother

\renewcommand		{\:}		{\colon}
\renewcommand		{\-}		{^{-1}}
\let\textoslash\o
\renewcommand		{\o}		{\circ}

\newcommand		{\adj}		{\dashv}

\ExplSyntaxOn
\NewDocumentEnvironment{adjunctions}{O{}}
{
	\cs_set_eq:cN {@arraycr} \farin_arraycr:
	\keys_set:nn { farin/adjunction } { #1 }
	\begin{array}
		{
			@{ \hspace { \dim_eval:n { \l_farin_left_shift_dim + \l_farin_padding_dim } } }
			r
			@{ {\farin_strut:} \l_farin_symbol_tl {} }
			l
			@{ \hspace { \dim_eval:n { \l_farin_right_shift_dim + \l_farin_padding_dim } } }
		}
	}
	{
	\end{array}
}
\keys_define:nn { farin/adjunction }
{
	leftshift       .dim_set:N = \l_farin_left_shift_dim,
	leftshift       .initial:n = 0pt,
	rightshift      .dim_set:N = \l_farin_right_shift_dim,
	rightshift      .initial:n = 0pt,
	padding         .dim_set:N = \l_farin_padding_dim,
	padding         .initial:n = 6pt,
	symbol          .tl_set:N  = \l_farin_symbol_tl,
	symbol          .initial:n = \longrightarrow,
	verticalspacing .dim_set:N  = \l_farin_vertspac_dim,
	verticalspacing .initial:n = {3pt},
}
\cs_new_protected:Npn \farin_strut:
{
	\vrule height \dim_eval:n { \ht\strutbox + 1.2\l_farin_vertspac_dim }
	depth  \dim_eval:n { \dp\strutbox + \l_farin_vertspac_dim }
	width 0pt
}
\makeatletter
\exp_args:NNo \cs_new:Npn \farin_arraycr:
{
	\@arraycr\hline
}
\makeatother
\ExplSyntaxOff

\renewcommand		{\.}		{\cdot}

\newcommand		{\x}		{\times}

\DeclareMathOperator*	{\otimesvariable}{%
			\mathchoice {\raisebox{.85pt}{$\displaystyle\otimes$}}
						{\raisebox{.85pt}{$\otimes$}}
						{\raisebox{0.7pt}{$\scriptstyle\otimes$}}
						{\raisebox{0.2pt}{$\scriptscriptstyle\otimes$}}
						}
\newcommand		{\tensor}	{\otimesvariable}

\newcommand		{\ox}		{\tensor}
\newcommand		{\ot}		{\direct}

\newcommand		{\Direct}	{\bigoplus}

\newcommand		{\ext}		{\exterior}
\newcommand		{\susp}		{\Sigma}

\newcommand		{\bul}		{\bullet}

\DeclareMathOperator	{\im}		{im }

\DeclareMathOperator	{\coker}	{coker }

\DeclareMathOperator	{\Tor}		{Tor}

\makeatletter
\newbox\xrat@below
\newbox\xrat@above
\newcommand		{\xrightarrowtail}[2][]	{%
						  \setbox\xrat@below=\hbox{\ensuremath{\scriptstyle #1}}%
						  \setbox\xrat@above=\hbox{\ensuremath{\scriptstyle #2}}%
				  \pgfmathsetlengthmacro{\xrat@len}{max(\wd\xrat@below,\wd\xrat@above)+.6em}%
  						\mathrel{\tikz [>->,baseline=-.55ex]
              					   \draw (0,0) -- node[below=-2pt] {\box\xrat@below}
                            					    node[above=-2pt] {\box\xrat@above}
                    						   (\xrat@len,0) ;}
						}
\newbox\xrat@below
\newbox\xrat@above
\renewcommand		{\xtwoheadrightarrow}[2][]{%
						  \setbox\xrat@below=\hbox{\ensuremath{\scriptstyle #1}}%
						  \setbox\xrat@above=\hbox{\ensuremath{\scriptstyle #2}}%
				  \pgfmathsetlengthmacro{\xrat@len}{max(\wd\xrat@below,\wd\xrat@above)+.6em}%
						 \mathrel{\tikz [->>,baseline=-.55ex]
					                 \draw (0,0) -- node[below=-2pt] {\box\xrat@below}
					                                node[above=-2pt] {\box\xrat@above}
						                       (\xrat@len,0) ;}
		       				}
\makeatother
\newcommand		{\xmono}	{\xrightarrowtail}

\newcommand		{\xepi}		{\xtwoheadrightarrow}
\newcommand		{\epi}		{\xepi{\phantom{\ \, }}}

\newcommand		{\longto} 	{\longrightarrow}
\newcommand		{\lt}		{\longto}
\newcommand		{\xtoo}		{\xrightarrow} 
\newcommand		{\from}		{\leftarrow}
\newcommand		{\longfrom}	{\longleftarrow}

\newcommand		{\lmt}		{\longmapsto}

\newcommand		{\simto}	{\xrightarrow{\sim}}

\newcommand		{\longsimto}	{\os\sim\longto}
\newcommand		{\isoto}	{\longsimto}

\newcommand		{\vertsim}	{\rotatebox{90}{$\sim$}}

\newcommand		{\hmt}		{\simeq}
\newcommand		{\iso}		{\cong}

\newcommand		{\Z}		{\bb Z}

\newcommand		{\C}		{\bb C}

\DeclareMathOperator	{\id}		{id}

\newcommand		{\CGA}		{\textsc{cga}\xspace}

\newcommand		{\DGA}		{\textsc{dga}\xspace}
\newcommand		{\DG}		{\textsc{dg}\xspace}
\newcommand		{\DGAs}		{\textsc{dga}s\xspace}
\newcommand		{\DGC}		{\textsc{dgc}\xspace}
\newcommand		{\DGCs}		{\textsc{dgc}s\xspace}
\newcommand		{\CDGA}		{\textsc{cdga}\xspace}
\newcommand		{\CDGAs}	{\textsc{cdga}s\xspace}
\newcommand		{\HGA}		{\textsc{hga}\xspace}
\newcommand		{\HGAs}		{\textsc{hga}s\xspace}

\newcommand		{\SHC}		{\textsc{shc}\xspace}
\newcommand		{\Ai}		{$A_\infty$}

\newcommand		{\kk}		{k}

\renewcommand 		{\H}		{H^*}

\newcommand		{\EMSS}		{Eilenberg--Moore spectral sequence\xspace}

\usepackage{tracklang}
\usepackage{multirow,accents}
\usepackage{array}\newcolumntype{R}{>{$}l<{$}}
\usepackage{scalerel}[2016/12/29]
\usepackage{chngcntr}
\counterwithin{table}{section}
\counterwithin{equation}{section}

\theoremstyle{definition}
\newtheorem*{layout*}{Outline}

\usepackage{silence}
\colorlet{jlabel}{Crimson!75!Black}		
\colorlet{jred}{Crimson!75!Black}	
\colorlet{nonisored}{Crimson!90!Black}	
\colorlet{jhmt}{Dandelion!80!Indigo}	
\colorlet{jcmp}{Black!35!White}
\colorlet{jblue}{RoyalBlue!80!Black}	
\colorlet{jviolet}{Violet!60!Black}
\colorlet{compromise}{rgb:Navy!85!Indigo,3;white,1}	
\colorlet{jQ2}{rgb:Indigo!90!Black,3;red,1}	
\colorlet{jA}{RoyalBlue!40!Indigo}		
\colorlet{jX3}{BlueViolet!100!Black}
\colorlet{jX2}{BlueViolet!100!Black}
\colorlet{jX2p}{Dandelion!50!Magenta}	
\colorlet{jgreen}{SeaGreen!85}
\colorlet{jcyan}{Cyan!70!Black}

\let		\epsilon	\varepsilon
\newcommand		{\B}		{\mathbf{B}}
\renewcommand	{\W}		{\BW}

\renewcommand	{\ot}		{^{\ox \mnn 2}}

\newcommand		{\htn}		{\vphantom{x^{x^{x^{x}}}}}
\newcommand		{\bariso}	{\iota}
\newcommand		{\predbariso}	{f}
\newcommand		{\WBbariso}	{\wt\bariso}
\newcommand		{\predsquism}	{g}
\newcommand		{\evquism}	{i}
\newcommand		{\Pb}		{P^\bul}

\renewcommand	{\susp}			{s}
\newcommand		{\desusp}		{s^{-1}}
\def				\iter#1#2{#1^{[#2]}}
\newcommand		{\muA}		{\mu_{\mn A}}
\renewcommand	{\C}		{C^*}
\newcommand		{\I}		{I^*}

\newcommand		{\Algs}		{\sans{DGA}}
\newcommand		{\Coalgs}	{\sans{DGC}}
\newcommand		{\Tw}		{\sans{T{\mn}w}}
\newcommand		{\Cochains}	{\sans{Ch}^\bul}
\newcommand		{\GM}		{\sans{Mod}}

\renewcommand	{\EMSS}		{\textsc{emss}\xspace}
\newcommand		{\SHCA}		{\SHC-algebra\xspace}

\newcommand		{\WHC}		{\textsc{whc}\xspace}
\newcommand		{\WHCA}		{\WHC-algebra\xspace}

\newcommand		{\TOR}		{\mathrm{TOR}}
\renewcommand{\ext}{\mr{ext}}
\newcommand		{\EXT}		{{\mr{ext.}}}

\DeclareMathOperator*{\T}{%
	\mathchoice {\raisebox{.85pt}{$\displaystyle\underline{{\otimes}}$}}
	{\mathbin{{\raisebox{.85pt}{$\underline{{\otimes}}$}}}}
	{\mathbin{{\raisebox{.7pt}{$\scriptstyle\underline{{\otimes}}$}}}}
	{\mathbin{{\raisebox{.2pt}{$\scriptscriptstyle\underline{{\otimes}}$}}}}
}

\DeclareFontFamily{U}{wncy}{}
\DeclareFontShape{U}{wncy}{m}{n}{<->wncyr10}{}
\DeclareSymbolFont{mcy}{U}{wncy}{m}{n}
\DeclareMathSymbol{\Sha}{\mathord}{mcy}{"58}

\def				\BW			{\boldsymbol{\Omega}}

\newcommand{\Tr}[1]{	\us{#1}	{\mathrm{Tor}}	}

\newcommand{\BA}{\B A}
\newcommand{\BX}{\B X}
\newcommand{\BY}{\B Y}

\newcommand{\tinystar}{{\scalebox{.5}{$*$}}}
\newcommand{\medstar}{{\raisebox{-.25ex}{\scalebox{.75}{$*$}}}}

\DeclareMathVersion{normal2}

\begin{document}

\title
	{%
	Products on Tor
	}
\author{Jeffrey D. Carlson}

\maketitle

\begin{abstract}
	In 1974 work establishing the collapse of
	certain Eilenberg--Moore spectral sequences,
	Munkholm constructs, in passing,
	a bilinear multiplication operation on 
	Tor of a triple of $A_\infty$-algebras.
	In 2020, the present author,
	pursuing a multiplicative collapse result extending Munkholm's,
	studied	a variant of this product,
	without actually showing it agrees with Munkholm's.
	In 2019, Franz had defined
	a weak product on the two-sided bar construction
	of a triple of $A_\infty$-algebras under similar hypotheses,
	with which this author proved a related collapse result,
    	but without investigating the properties of the induced
	product on Tor.

	The present work 
	demonstrates that the two products
	on Tor agree 
	and are induced by the product of Franz.
\end{abstract}

At the beginning of homological algebra lie the 
derived functors $\Tor^i$ of the tensor product.
When $M \from A \to N$ are maps of commutative graded algebras
(\defm{\CGA}s),
the graded groups $\Tor^i_A(M,N)$ fit into a bigraded ring
$\defm{\Tor_A(M,N)}$
because the multiplications
$A \ox A \lt A$ and so on are themselves ring maps. 
When $M \from A \to N$ are maps of noncommutative 
differential graded algebras (\defm{\DGA}s),
there is still an appropriate notion of proper projective resolution
for differential graded $A$-modules,
and accordingly a tensor product whose 
derived functors are again written $\Tor^i_A(M,N)$,
reducing to the more classical notion when the differentials are zero.
The definition from before typically no longer yields
a ring structure on Tor,
but when the input rings are (noncommutative) cochain algebras 
$\C(X) \from \C(B) \to \C(E)$,
there is a different  ring structure
arising from the Eilenberg--Zilber theorem, 
which can be seen as a homotopy-commutativity
property of (co)chains.

Assuming a weaker homotopy-commutativity hypothesis
on the input \DGAs,
Hans J{\textoslash}rgen Munkholm 
defined a product on Tor generalizing these two examples, 
and promptly abandoned it~\cite[\SS9]{munkholm1974emss}.
In recent work%
    ~\cite{carlson2022munkholm}, the present author revived this product,
to various ends that need not detain us here.\footnote{\ 
    In particular, the original content of this paper
    is independent of that of the previous works
    apart from background and definitions. 
}
The definition of the revived product is actually a simplification
of Munkholm's original, 
and owing to length considerations,
proof that the definitions are equivalent
is not actually included in that work,
it having seemed more important at the time 
to establish the relevant properties of \emph{some} product.
Similarly, previous work of the author 
employed a product due to Franz
on a two-sided bar construction $\B(X,A,Y)$,
defined under similar hypotheses%
~\cite{carlsonfranzlong},
which also induces a product on $\Tor_A(X,Y)$
under mild flatness hypotheses,
whose properties went unexplored in favor of applications.

In the present paper, we show 
that the product of our previous work
is indeed Munkholm's (\Cref{thm:Jeff-equiv-Munkholm}).
This requires us to recognize 
a connection 
that previously had gone unremarked
between three classically defined
natural transformations of \DGAs
(\Cref{thm:WDg}).\footnote{\ 
    This result may be of independent interest
    because it can be used to prove the crucial Proposition IV.6.1 
    of Husemoller--Stasheff--Moore's independent collapse 
    paper~\cite[p.~179]{husemollermoorestasheff1974}, 
    which in the original writing is supported 
    by a statement,
    Proposition IV.5.7,
    that is unfortunately not true.
}
We then show that
so long as $\B(X,A,Y)$ does compute Tor,
Franz's product on $\B(X,A,Y)$
induces Munkholm's product on $\Tor_A(X,Y)$ 
(\Cref{thm:Jeff-equiv-Franz}).

\section{Algebras, coalgebras, and twisting cochains}\label{sec:CDGA}
Prerequisites are as in the predecessor~\cite{carlson2022munkholm},
but we run rapidly over some highlights.

\begin{notation}
We take tensors and Homs over a fixed commutative base ring $\kk$ with unity
and consider \emph{\textbf{nonnegatively-graded}}
cochain complexes $(C,\defm d)$,
with differentials $d$ \textbf{\emph{increasing degree}} by $1$,
writing $\defm\Cochains$ for the category thereof.
The Koszul sign convention is always in force.
We use only augmented
differential nonnegatively-graded $\kk$-algebras (henceforth \textcolor{RoyalBlue}{\DGA}s)
$(A,+,0,{d_A},\defm{\muA},\defm{\h_A},\defm{\e_A})$
with augmentation ideal $\ker \e_A = \defm{\ol A} \iso \coker \h_A$
and coaugmented, cocomplete, differential nonnegatively-graded $\kk$-coalgebras (\textcolor{RoyalBlue}{\DGC}s)
$(C,+,0,d_C,\defm{\D_A},\defm{\e_A},\defm{\h_A})$
with coaugmentation coideal $\coker \h_A = \defm{\ol C} \iso \ker \e_A$.
The notions of homomorphism are as expected
and the corresponding categories are written 
$\defm{\Algs}$ and $\defm{\Coalgs}$ respectively.
The base ring $\kk$ itself lies in both.
A commutative \DGA is a \textcolor{RoyalBlue}{\CDGA}.
We use the terms \emph{\DG ($\kk$-)module} and \emph{cochain complex}
interchangeably.
\end{notation}

\begin{definition}\label{def:cup}
Given two graded $\kk$-modules $C$ and $A$,
we denote by $\defm{\GM_n}(C,A)$ 
the $\kk$-module of $\kk$-linear maps $f$ sending each $C_j$ to $A_{j+n}$,
and set the degree $\defm{|f|}$ to $n$ for such a map.
The hom-set $\GM(C,A) = \Direct_{n \in \Z} \GM_n(C,A)$ 
then becomes itself a graded $\kk$-module,
an internal Hom in the category $\defm\GM$ of graded $\kk$-modules.
If $C$ and $A$ are cochain complexes,
then $\GM(C,A)$ becomes
a cochain complex under the
differential~$D = {d_{\GM(C,A)}}$
given by $\defm{Df} \ceq d\mn_A f - (-1)^{|f|}f\mnn d_C$~{\cite[\SS1.1]{munkholm1974emss}}.
The cochain maps $C \lt A$ are precisely $\ker D$.
If $C$ is a \DGC and $A$ a \DGA,
then $\GM(C,A)$ becomes a \DGA
under the \defd{cup product}
$\defm{f \cup g} \ceq \muA(f \tensor g)\D_C$%
~{\cite[\SS1.8]{munkholm1974emss}},
with unity $\defm{*} \ceq \h_A \e_C$.
An element  $t \in \GM_1(C,A)$
satisfying the three conditions
\[
\e_{\mn A} t = 0 = t \h_C,\qquad Dt = t \cup t
\]
is called a \defd{twisting cochain}%
~\cite[\SS{1.8}]{husemollermoorestasheff1974}%
\cite[Prop.~3.5(1)]{husemollermoorestasheff1974}%
\cite[\SS\SS1.5,\,4]{proute2011};
we write $\defm{\Tw(C,A)}$ for the additive group of these.
Given
$(g,t,f) \in \Coalgs(C',C) \x \Tw(C,A) \x \Algs(A,A')$,
the maps $tg$, $f\mnn t$, and $f\mnn tg$ are again twisting cochains.
\end{definition}

\begin{definition}\label{def:bar}
For each \DGA $A$,
there is a \emph{final} twisting cochain
$\defm{t^A}\: \defm {\B A} \lt A$
defined by the property that any 
twisting cochain $t\: C \lt A$
factors uniquely through a \DGC
map $\defm{g_t}\: C \lt \B A$
such that $t = t^A \o g_t$. 
Here the cocomplete \CGA 
$\B A$ is the familiar (normalized) \defd{bar construction},
which gives the object component of a functor $\defm \B\:\Algs \lt \Coalgs$%
~{\cite[\SS1.6]{munkholm1974emss}\cite[\SS2.5]{proute2011}}.
Write $\defm {U_\Algs}\: \Algs \lt \Cochains$
and $\defm {U_\Coalgs}\: \Coalgs \lt \Cochains$
for the forgetful functors.
The tautological twisting cochain $t^{(-)}\: U_\Coalgs \o \B \lt U_\Algs$
is a natural transformation of functors $\Algs \lt \Cochains$.
We denote this conversion in the input--output ``deduction rule'' format:
\[
\begin{adjunctions}
	g_t\: C & \B A\\
	t\: C & A \mathrlap.
\end{adjunctions}
\]
\nd Explicitly, 
$\BA$ is the tensor coalgebra on 
the desuspension $\defm{\desusp \ol A}$ of $\ol A$,
equipped with the sum of the tensor differential
and the unique coderivation extending
the ``bar-deletion'' map.
\end{definition}

\begin{observation}\label{thm:BA-structure}
Writing $\B_n A = (\desusp\ol A)^{\otimes n}$
for the summands of $\B A$,
and $\defm{\ol \D} b = \D b - 1 \ox b - b \ox 1$
for the reduced comultiplication,
$\bar\Delta \B_n A$ 
lies in the sum of $\B_p A \ox \B_q A$
for $p+q = n$ and $p,q \geq 1$.
The kernels  of
$\iter {\bar \D{}}{n}$ are
$\smash{\Direct_{p < n} \B_p A}$,
which form a filtration by subcomplexes;
in particular, $\B_1 A \iso \desusp \ol A$ is a subcomplex.
The tautological twisting cochain $t^A\: \B A \lt A$
factors through a cochain isomorphism 
$\B_1 A \isoto \ol A$ of degree $1$,
whose inverse is $\defm\desusp\:  \ol A \simto \desusp \ol A = \B_1 A$.
\end{observation}

\bdefn
For each a \DGC $C$, there is a twisting cochain
$\defm{t_C}\: C \lt \defm {\W C}$
\emph{initial} in the sense that any twisting cochain 
$t\: C \lt A$
factors uniquely through a $\DGA$ map 
$\defm{f^t}\: \W C \lt A$ 
such that $t = f^t t_C$.
The \DGA $\W C$ is referred to as the \defd{cobar construction},
and gives the object component of a functor $\defm\W\:\Coalgs \lt \Algs$%
~{\cite[\SS1.7]{munkholm1974emss}}.
Its underlying algebra is 
the tensor algebra $\Direct_{n \geq 0} \W_n C$,
where $\defm{\W_n C} = (\susp \ol C)^{\otimes n}$.
Write $\defm {V_\Algs}\: \Algs \lt \GM$
and $\defm {V_\Coalgs}\: \Coalgs \lt \GM$
for the forgetful functors.
The tautological twisting cochain $t_{(-)}\: V_\Coalgs \lt V_\Algs \o \W$
is a natural transformation of functors $\Coalgs \lt \Cochains$.

The two functors $\W \adj \B$ form an adjoint pair~{\cite[\SS1.9--10]{munkholm1974emss}}.
We will have frequent recourse to the unit and counit of the adjunction $\W \adj \B$,
\[
	\defm\h\: \id \lt \B\W
\qquad\qquad\mbox{and}\qquad\qquad
	\defm\e\: \W\B \lt \id
\]
respectively. 
These are both natural quasi-isomorphisms 
and homotopy equivalences on the level of \DG modules%
~\cite[Thm.~II.4.4--5]{husemollermoorestasheff1974}%
\cite[Cor.~2.15]{munkholm1974emss}%
\cite[Lem.~1.3.2.3]{LH}.
\edefn


The adjunction interacts with the tautological twisting cochains as follows.

%
\begin{lemma}\label{thm:twisting-adjunction}
	For a \DGA $A$ and a \DGC $C$, 
	one has $\e \o t_{\B A} = t^A\: \BA \lt A$
	and $t^{\W C} \o \h = t_C\: C \lt \W C$.
\end{lemma}

\section{The tensor product}\label{sec:tensor}

The functor $\B\: \Algs \lt \Coalgs$ is lax monoidal with respect to the 
monoidal structure given on both categories by the appropriate tensor products,
and $\W\: \Coalgs \lt \Algs$ is lax comonoidal.

\begin{definition}[See Husemoller \emph{et al.}~{\cite[Def.~IV.5.3]{husemollermoorestasheff1974}}]%
\label{def:shuffle}
There exist natural transformations 
\[
		\defm\nabla \: \B A_1 \ox \B A_2 \lt \B(A_1 \ox A_2),	
\qquad\qquad
		\defm\gamma \: \W( C_1 \ox C_2) \lt \W C_1 \ox \W C_2
\]
\nd of functors $\Algs \x \Algs \lt \Coalgs$ and $\Coalgs \x \Coalgs \lt \Algs$,
respectively, the 
\defd{shuffle maps}, 
determined by the twisting cochains
%
\[
			t^{A_1 \ox A_2}\nabla  = t^{A_1} \ox \h_{A_2}\e_{\B A_2} +
			\h_{A_1}\e_{\B A_1} \ox t^{A_2},
\qquad\qquad
		\gamma t_{C_1 \ox C_2}  = t_{C_1} \ox \h_{\W C_2}\e_{C_2} +
				\h_{\W C_1}\e_{C_1} \ox t_{C_2}\mathrlap.
\]
These are homotopy equivalences of cochain complexes and hence quasi-isomorphisms.
\end{definition}

Written out in terms of bar-words,
$\nabla(b_1 \ox b_2)$ is a sum of shuffle permutations 
of the letters of $b_1$ and $b_2$,
so values of $\nabla$ exhibit symmetry
with respect to shuffles of tensor-factors.

\begin{observation}\label{lem:symmetry}
Let $A_1$ and $A_2$ be \DGAs.
For $b_1 \ox b_2 \in \B_{m} A_1 \ox \B_{n} A_2$,
the $\nabla\ot$-image of the summands
of $\D_{\B A_1 \ox \B A_2}(b_1 \ox b_2)$
lying in $\B_{m} A_1 \ox \kk \ox \kk \ox \B_{n} A_2$
and $\kk \ox \B_{n} A_2 \ox \B_{m} A_1 \ox \kk$, 
respectively
$(b_1 \ox 1) \otimes (1 \ox b_2)$ and 
$(-1)^{|b_1||b_2|}(1 \ox b_2) \otimes (b_1 \ox 1)$,
is
\[
\nabla(b_1 \ox 1) \ox \nabla(1 \ox b_2) + (-1)^{|b_1||b_2|}
	\nabla(1 \ox b_2) \ox \nabla(b_1 \ox 1)\mathrlap.
\]
\end{observation}


We will require some more terminology 
to introduce the important natural transformation $\psi$ of \Cref{def:psi},
the details of whose construction play a role
in the proof of \Cref{thm:WDg}.

\begin{definition}[{\cite[\SS2.1]{munkholm1974emss}}]\label{def:TEX}
	A \defd{trivialized extension}\footnote{\ 
		We follow Munkholm in this usage.
		\emph{Contraction},
		\emph{(strong) homotopy retract datum},
		and \emph{SDR-data} are all common 
		in the literature 
		when $\wt A$ and $A$ are merely assumed
		\DG modules.
	} 
	is an assemblage of maps 
	$\smash{
	\xymatrix@R=1.5em{
		\wt A 	
		\ar@<-0.3ex>[r]_{p} 
		\ar@(dl,ul)[]^{h} 
		&
		A 		
		\ar@<-0.3ex>[l]_{i} 
		}
		}
	$
with $p$ a \DGA map, 
	$i$ a degree-$0$ \DG module section,
	and $h \in \GM_{-1}(\wt A,\wt A)$
	a cochain homotopy satisfying $Dh = \id{} - i\mnn p$
	and such that, moreover,
	the compositions $ph$, $hh$, $hi$ 
	vanish.\footnote{\ 
		That $h^2 = 0$ 
		actually follows from the other equations.
		}
\end{definition}

The homotopy $h$ allows us to promote $i$
to a \DGC map $\B A \lt \B \wt A$.

\begin{lemma}[Homotopy transfer theorem for \DGAs~{\cite[Prop.~2.2]{munkholm1974emss}}]\label{thm:promotion}
	Let a trivialized extension be given as in \Cref{def:TEX}.
	Then there exists a twisting cochain $\defm{t^i}\: \B A \lt \wt A$
	such that $pt^i = t^A\: \B A \lt A$.
This $t^i$, given recursively by $t^i = h(t^i \cup t^i) + it^A$,
then induces a \DGC map 
	$g_{t^i}\:\B A \lt \B \wt A$
	and a \DGA map 
	$f^{t^i} \: \W\B A \lt \wt A$.%
	\footnote{\
		Munkholm does not write out a proof, but it is not completely trivial.
	\begin{proof}
	We first explain the recursive prescription for $t = t^i$.
	By \Cref{thm:BA-structure}, 
	we may recursively define $t_n = t|_{\B_n A}$.
	Our formula sets $t_0 = 0$
	and $t_1 = it^A|_{\B_1 A}$,
	and for $n \geq 2$ takes $t_n(b) = h(t_{< n} \cup t_{< n})(b)$,
	which avoids circularity because 
	$t \cup t$ annihilates $1 \ox b$ and $b \ox 1$
	by the ``$t_0 = 0$'' clause
	and $\bar\D_{\B A}$ takes $\B_n A$ to 
	$\B_{< n} A \ox \B_{< n} A$,
	where $t \cup t$ is already defined. 
	Evidently $pt_n = 0$ for $n > 1$ since $ph = 0$,
	and $pt_0 = p0 = 0$, while $pt_1 = pit^A = t^A$.
	
	The proof $t$ is a twisting cochain is by induction.
	That $Dt_0 = 0 = (t \cup t)|_{\B_0 A}$ is trivial.
	We have $Dt_1 = dit^A + it^A d = 0$
	because $i$ is a chain map of degree $0$ and $t^A|_{\B_1 A}$
	a chain map of degree $1$,
	and $t \cup t$ vanishes on $\B_1 A$ since $t_0 = 0$ and $\bar \D$ 
	vanishes on $\B_1 A \iso \ker \ol \D \,\mnn/\mn \im \eta_{\B A}$. 
	
	The case $n = 2$ is the interesting case.
	On the one hand one has $(t \cup t)|_{\B_2 A} = t_1 \cup t_1$.
	On the other, since $t^A$ is a twisting cochain,
	one has $dt^A + t^A d|_{\B_2 A} = t^A \cup t^A$,
	and using in order these facts,
	that $Dt_1 = 0$
	and $D(t_1 \cup t_1) = 0$,
	that $dh + hd = {\id} - i\mnn p$, 
	and that $pi = \id\mn_A$, we find
	\eqn{
	(Dt)|_{\B_2 A} 	&= dt_2 + t_{\leq 2} d \\
				&= dh(t_1 \cup t_1) + dit^A + h(t_1 \cup t_1)d + it^A d\\
				&= dh(t_1 \cup t_1)
				+ hd(t_1 \cup t_1)
				+ idt^A + it^A d
				\\
				&= (dh+hd)(t_1 \cup t_1) + i(t^A \cup t^A) \\
				&= t_1 \cup t_1\, - \,i\mnn p(it^A \cup it^A)\, + \,it^A \cup it^A\\
				&= t_1 \cup t_1\mathrlap.
	}
	
	For $n \geq 3$ the proof is tautological:	
	we have $t_n = h(t \cup t)$
	since $d \B_n A$ lies in $\B_{\geq n-1} A$
	and hence is annihilated by $it^A$,
	and so $\im (t \cup t)|_{\B_n A}$ lies in the ideal generated by $\im h$,
	which $p$ annihilates, hence inductively 
	\eqn{
		Dt   	&= dh(t \cup t) + h(t \cup t)d			\\
				&= (\id{} - i\mnn p - hd)(t \cup t) + h(t \cup t)d	\\
				&= t \cup t\, - \,h \, D(t \cup t)				\\
				&= t \cup t\, -\, h\, D\mnn\big(\mspace{-1.5mu} D(t)\mspace{-1.5mu}\big) \\
				&= t \cup t\mathrlap.
		\qedhere
	}
	\end{proof}
	}
\end{lemma}

We will not make real use of 
morphisms of trivialized extensions---%
pairs of \DGA maps making the 
expected three squares for $p$, $i$, $h$ 
commute~\cite[\SS2.1]{munkholm1974emss}---%
but we will need one key example.

\bex[The universal example~{\cite[Prop.~2.14]{munkholm1974emss}}]%
\label{def:TEX-WBA}
Given a \DGA $A$,
there is a unique section $\defm{i_A}\: A \lt \W\B A$ of 
$\e\: \W\B A \lt A$
defined to be unital and to restrict to $t_{\B A}\o \desusp$ on $\ol A$.
Along with a certain homotopy $h$ in $\GM_{-1}$$(\W\B A,\W\B A)$ 
we will not be explicit about, 
$\e$ and $i_A$
can be shown to give a trivialized extension.
The only detail we will need about $h$
is that it is constructed inductively using
an decomposition
$\Direct_{j=0} S_j$ of $\W\BA$
by graded submodules $\defm{S_j}$
such that $S_0$ is annihilated by $h$ 
and contains $\W_0\BA$ and $\W_1 \BA$
(cobar-words of length $0$ and $1$)%
~\cite[p.~17]{munkholm1974emss}.

The cochain $t^{i_A} = h(t \cup t) + i_At^A\: \B A \lt \W\B A$
defined recursively from $t^A\: \B A \lt A$
as in \Cref{thm:promotion}
works out to be the tautological $t_{\B A}$.
Given another trivialized extension 
$p\:\wt A \to A$ with section $i\: A \to \wt A$,
by \Cref{thm:promotion} 
there is an induced \DGA map 
$\smash{f^{t^{i}}}\: \W\B A \lt \wt A$ 
satisfying $p\mn f^{t^i}= \e\: \W\B A \lt A$.
\eex

\begin{theorem}[{%
\cite[Prop.~IV.5.5]{husemollermoorestasheff1974}%
\cite[{$k_{A_1,A_2}$, p.~21, via Prop.~2.14}]{munkholm1974emss}%
}]%
\label{def:psi}
	There exists a natural transformation 
		\[\defm \psi\: \W\B(A_1 \ox A_2) \lt \W\B A_1 \ox \W\B A_2\]
	of functors $\Algs \x \Algs \lt \Algs$.
	This transformation satisfies 
	\[
		(\e_{\mn A_1} \ox \e_{\mn A_2}) \o \psi 
			= 
		\e_{\mn A_1 \ox A_2}\: \W\B(A_1 \ox A_2) \lt A_1 \ox A_2		
	\]
	and reduces to the identity if $A_1$ or $A_2$ is $\kk$.%
    \footnote{\ 
	Husemoller--Moore--Stasheff's version of $\psi$
	is defined using a splitting result (IV.2.5)
    depending on certain objects
    being injective in an appropriate sense.
    It has the same categorical properties as Munkholm's map,
	but because we will need a somewhat more explicit 
    cochain-level description of $\psi$ in \Cref{thm:WDg},
     we employ Munkholm's formulation.
     }
\end{theorem}
\begin{proof}
	Granting the claims of \Cref{def:TEX-WBA},
	it is easy to check the data
	\eqn{
	\e_{A_1} \ox \e_{A_2}\: 
		&\W\B A_1 \ox \W \B A_2 \lt A_1 \ox A_2\mathrlap,\\ 
	i_{1} \ox i_{2}\: 
		&A_1 \ox A_2 \lt \W\B A_1 \ox \W \B A_2\mathrlap,\\
\defm{h_\psi} \ceq	h_{1} \ox {\id} + i_{1}\e_{A_1} \ox h_{2}\:
		& \W\B A_1 \ox \W \B A_2 \lt \W\B A_1 \ox \W \B A_2
	}
	give a trivialized extension,\footnote{\ 
		More generally,
		Munkholm defines a tensor product
		of trivialized extensions~\cite[Prop.~2.10]{munkholm1974emss}.
	}
	so 
	\Cref{thm:promotion} 
	yields the required 
	\DGA map $\psi = 
	\smash{
		f^{\,\mn 
			t^{\,\mnn 
				i_{\mn\smash{1}}\! 
				\ox 
				i_{\mn\smash{2}}
			}
		} 
	}$,
	with associated twisting cochain 
	$\defm{t_\psi}\: \B(A_1 \ox  A_2) \lt \W\B A_1 \ox \W\B A_2$.
	Naturality follows from the naturality of 
	$\e$, $i$, and $h$.
	
	If $A_2 = \kk$, 
	then we may make the identifications
	$\W\B(A_1 \ox \kk) = \W \B A_1$
	and $\W \B A_1 \ox \W\B A_2 = \W\B A_1 \ox \kk = \W\B A_1$,
	so that $\e_{A_2}$ and $i_{2}$ are identified 
	with $\id_\kk$
	and $h_{2} i_{2} = 0$ forces $h_{2} = 0$,
	and make the identification $h_\psi = h_{1}$.
	Thus this trivialized extension reduces
	to the initial example $\W\B A_1 \lt A_1$
	of \Cref{def:TEX-WBA}.
	Now $\psi\: \W\B A_1 \lt \W\B A_1$ 
	is induced from 
	$\smash{t_\psi}\: \B A_1 \lt \W\B A_1$,
	recursively defined
	by $t_\psi = h_{1}(t_\psi \cup t_\psi) + i_{1}t^{A_1}$,
	but then $t_\psi$ agrees
	with $t^{i_{1}}$ from the universal
	\Cref{def:TEX-WBA},
	which we have stated is 
	the tautological twisting cochain $t_{\B A_1}$
	whose associated \DGA map is 
	$f^{t_{\B A_1}} = \id_{\W \B A_1}$.
	The proof if instead $A_1 = \kk$ is symmetrical.
\end{proof}


%

We will frequently consider \DGC maps 
$\B A \lt \B B$
between bar constructions on \DGAs,
which can be seen as a sort of up-to-homotopy version
of \DGA maps from $A$ to $B$,
  sometimes called \emph{\Ai-algebra maps}.
The natural transformation $\psi$ allows us to take tensor products of such maps.

\begin{definition}[{\cite[Prop.~3.3]{munkholm1974emss}}]\label{def:internal-tensor}
	Let $A_1,A_2,B_1,B_2$ be \DGAs and 
	$g_j\: \B A_j \lt \B B_j$ be \DGC maps for $j \in \{1,2\}$.
	Then we define the \defd{internal tensor product} 
	$\defm{g_1 \T g_2}\: \B(A_1 \ox A_2) \lt \B(B_1 \ox B_2)$ 
	by 
	\[
	\B(A_1 \ox A_2) 
		\xtoo{\!\h\!} 
	\B\W\B(A_1 \ox A_2)
		\xtoo{\B\psi}
	\B(\W\B A_1 \ox \W\B A_2)
		\xtoo{\B(\e \,\W g_1 \ox \,\e \, \W g_2)}
	\B(B_1 \ox B_2)
		\mathrlap.
	\]
\end{definition}

We will use one relation between the counit $\e$
and the internal tensor product.

\begin{lemma}[{\cite[p.~49, top]{munkholm1974emss}}]\label{thm:epsilon-T}
	Let $A_j$ and $B_j$ be \DGAs and $g_j\: \B A_j \lt \B B_j$ be \DGC maps
	for $j \in \{1,2\}$.
	Then one has
	\[
		\e\o\W(g_1 \T g_2) = (\e \ox \e)\o(\W g_1 \ox \W g_2)\o\psi
		\:
		\W\B(A_1 \ox A_2) \lt B_1 \ox B_2
		\mathrlap.
	\]
\end{lemma}
\begin{proof}
	We chase a commutative diagram.
\[
\resizebox{6.5in}{!}{
\xymatrix@C=2.75em@R=4em{
	\W\B(A_1 \ox A_2) \ar[d]_{\W\h} \ar[r]^(.475){\psi}&
	\W\B A_1 \ox \W\B A_2 \ar[rr]^(.46){\W g_1 \ox \W g_2}&&
	\W\B B_1 \ox \W\B B_2 \ar[r]^(.575){\e \ox \e} &
	B_1 \ox B_2	\\
	\W\B\W\B(A_1 \ox A_2) \ar[r]_(.475){\W\B\psi} &
	\W\B(\W\B A_1 \ox \W\B A_2) \ar[rr]_{\W\B(\W g_1 \ox \W g_2)} \ar[u]^\e&&
	\W\B(\W\B B_1 \ox \W\B B_2) \ar[r]_(.575){\W\B(\e \ox \e)} \ar[u]_\e&
	\W\B(B_1 \ox B_2) \ar[u]_\e
}
}
\]	
The composition along the top is the right-hand side of the display
and the composition along the bottom is the left-hand side, by \Cref{def:internal-tensor}.
The left square commutes since 
$\e \o \W\B\psi \o \W\h = \psi \o \e \o \W\h = \psi$
by naturality of $\e$ and the unit--counit identities 
    for the adjunction $\W \adj \B$,
and the other two squares commute by naturality of $\e$.
\end{proof}
%
%
%
%
%
%

\section{Homotopies and maps on Tor}\label{sec:Tor}

In this section we discuss notions of homotopy,
how to represent homotopies by maps into path objects,
and how to use such maps to define maps on Tor.

\begin{definition}%
	[{\cite[\SS1.11]{munkholm1974emss}\cite[\SS4.1]{munkholm1978dga}}]%
	\label{def:homotopy}
\ \\
	A \defd{homotopy}
	of \DGC maps $g_0,g_1\: C' \lt C$
	is a degree-$(-1)$ $\kk$-linear map $j\: C' \lt C$
	such that 
	\[
	\e_{C} j = 0,\qquad\quad
	j\eta_{C'} = 0,\qquad\quad
	d(j) = g_1 - g_0,\qquad\quad
	\D_{C} j = (g_0 \ox j + j \ox g_1)\D_{C'}
	\mathrlap.
	\]
	A \defd{homotopy}
	of twisting cochains $t_0,t_1\: C \lt A$
	is a degree-$0$ $\kk$-linear map $x\: C \lt A$ 
	such that 
	\[
	\e_{\mn A} x = \e_C,\qquad\quad
	x\h_{A} = \h_C,\qquad\quad
	\phantom{d(h) = f_0 - f_1,}\qquad
	d(x) = t_0 \cup x - x \cup t_1
	\mathrlap.
	\]
	A 
    \defd{homotopy}
	of \DGA maps $f_0,f_1\: A \lt A'$
	is a degree-$(-1)$ $\kk$-linear map $h\: A \lt A'$ 
	such that 
	\[
    \hphantom{{}^7}
    \e_{\mn A'} h = 0,\qquad\quad
	h\eta_{A} = 0,
    \qquad\quad
	d(h) = f_0 - f_1,\qquad\quad
	h\mu_{\mn A} = \mu_{\mn A'}(f_0 \ox h + h \ox f_1)
.%
    \footnote{\ %
		In the definition from our main source~\cite{munkholm1974emss}, 
		the unit and counit conditions are omitted;
		in later work dealing more specifically
		with $\Algs$ as a category,
		he includes them~\cite[4.1]{munkholm1978dga}.
		These details do not affect the definition(s)
        of the product here
        but are critical for the adjunction 
		to preserve homotopy and 
        for the path object
        to do what we require it to in the predecessor~\cite{carlson2022munkholm}.
	}
	\]
\end{definition}

\bs

These three notions 
compose as expected with maps in the appropriate categories
and are also interconverted by the adjunctions.

%
%
\begin{lemma}[{\cite[\SS1.11; Thm.~5.4, pf.]{munkholm1974emss}}]\label{thm:homotopy-adjunction}
	Suppose given a \DGC $C$ and a \DGA $A$.
	Then there are bijections of homotopies of maps
	\quation{\label{eq:homotopy-adjunction}
		\begin{aligned}
			\begin{adjunctions}
				\W C 	& A\\
				C 		& A\\
				C 		& \B A
			\end{adjunctions}
		\end{aligned}
	}
	The adjoint functors $\Coalgs$$\:\W \adj \B :\!\mn\Algs$ 
	also preserve the relation of homotopy.\footnote{\ 
		That the relation of homotopy is preserved
		is not to say that,
		for instance, if $j\: g_0 \hmt g_1\: C \lt C'$
		is a \DGC homotopy, then $\W j$ is a \DGA homotopy
		from $\W g_0$ to $\W g_1$,
		but that there exists a certain \DGA homotopy.
		A longer unpublished draft~%
  \cite{carlsonmunkholmlong}
  goes into detail about what this homotopy is;
    the primary sources do not seem to.
	}
\end{lemma} 

%
%
%

It is well known that the data of 
a homotopy $j\: g_0 \hmt g_1\: C \lt C'$ of maps of chain complexes (resp. \DGCs)
can be realized as single map $C \ox I \lt C'$,
where $\defm I$ is the complex 
$\kk\{u_{[0,1]}\} \to \kk\{u_{[0]},u_{[1]}\}$
of nondegenerate chains
in the standard simplicial structure on the interval $[0,1]$ (with 
the expected coproduct).
Munkholm~\cite[Thm.~5.4, pf.]{munkholm1974emss}
observed 
that the dual \DGA $\defm\I$ of normalized simplicial cochains on $[0,1]$
has the following dual property.

\begin{observation}\label{def:P} 
    The \DGA $\I$
is a free $\kk$-module of rank $3$ on a basis $v_0,v_1,e$
with $|v_0| = 1 = |v_1|$ and $|e| = 1$
and unity $1_{\I} = v_0 + v_1$.
The projections 
$
\defm{\pi_j} \: \I \ox A \epi \kk\{v_j\} \ox A \simto A
$
are \DGA quasi-isomorphisms with the property that
a \DGA homotopy $h\: f_0 \hmt f_1\: A' \lt A$
can be represented as a \DGA map
$
		\defm{h^P}\: A' \lt \I \ox A
		$
such that $\pi_j \o h^P = f_j$.

For compatibility with the predecessor~\cite{carlson2022munkholm},
we replace $\I \ox A$ with its 
quasi-isomorphic subalgebra
$\defm {PA} \ceq \kk\{1_{\I} \ox 1_{A}\} \,\oplus\, \I \ox \ol A$
in order that the path object be augmented in such a way that
(the restrictions of) the maps $\pi_j$ become
augmentation-preserving.\footnote{\ 
	To define the product,
	and make comparisons with other products,
	we could equally well use $\I \ox A$
	in this paper,
	but in order to study its properties 
        (which for instance requires forming $\B PA$),
	we do need $PA$.
}
\end{observation}

\counterwithin{figure}{section}
\counterwithin{definition}{section}
\counterwithin{theorem}{section}
\counterwithin{lemma}{section}
\counterwithin{corollary}{section}
\counterwithin{proposition}{section}
\counterwithin{equation}{section}
\counterwithin{notation}{section}
\numberwithin{remark}{section}

To use the path object, we recall the \emph{algebraic \EMSS},
a spectral sequence of K\"unneth type 
converging to differential Tor of a triple of \DGAs
and functorial in all three variables~\cite[XI.3.2]{maclane}.

\begin{lemma}[{\cite[Cor.~1.8]{gugenheimmay}%
		\cite[Theorem~5.4]{munkholm1974emss}}]
	\label{thm:Tor-quism}
Given a commutative diagram of \DGA maps
\quation{\label{eq:Tor-DGA-functoriality-squares}
	\begin{aligned}
		\xymatrix@C=2em{
			X' \ar[d]_(.45)u
			&	\ar[l]_{\phi_{X'}} \ar[r]^{\phi_{Y'}}
			A'
			\ar[d]|(.375)\hole|(.45)f|(.525)\hole
			& 	Y\mathrlap'
			\ar[d]^(.45)v\\
			X
			& A
			\ar[l]^{\phi_X} \ar[r]_{\phi_Y}	
			& Y\mathrlap,
		}
	\end{aligned}
	}
	 there is induced a map of algebraic {\EMSS}s
	from that of $(X',A',Y')$ to that of $(X,A,Y)$,
	converging
	to the functorial map
	$
	\Tor_f(u,v)\:	\Tor_{A'}(X',Y') \lt \Tor_{A}(X,Y)
	$
	of graded modules.
	Hence when the maps $f$, $u$, $v$
	are quasi-isomorphisms,
    $\Tor_f(u,v)$ is an isomorphism.
\end{lemma}
%
We will need to expand the notion of a map of Tors.
 
\begin{lemma}[{\cite[Thm.~5.4]{munkholm1974emss}}]%
	\label{thm:homotopy-Tor-map}
	Let \DGA maps as in \eqref{eq:Tor-DGA-functoriality-squares}
	be given such that the squares commute up to \DGA homotopies
	$h_X\: u \o \phi_{X'} \hmt \phi_X \o f$ and
	$h_Y\: v \o \phi_{Y'} \hmt \phi_Y \o f$.
	Then there is induced a map 
	\[
	\defm{\Tor_{f}(u,v;h_X,h_Y)}\:
	\Tor_{A'}(X',Y')
	\lt
	\Tor_A(X,Y)
	\]
	of graded modules
	which is a quasi-isomorphism if each of $u$, $f$, and $v$ is.
\end{lemma}
\begin{proof}
	Letting $h_X^P\: A \lt PX'$
	and $h_Y^P\: A \lt PY'$
	be the \DGA representatives
	for the homotopies $h_X$, $h_Y$
	 described in \Cref{def:P},
	the following diagram commutes by definition:
	\quation{\label{eq:Tor-DGA-homotopy-squares}
		\begin{aligned}
			\xymatrix@C=1.90em{
				X' \ar[d]_(.45)u
				&	\ar[l]_{\phi_{X'}} \ar[r]^{\phi_{Y'}}
				A'
				\ar@{=}[d]
				& 	Y\mathrlap'
				\ar[d]^(.45)v\\
				X
				& A'\ar@{=}[d]
				\ar[l] \ar[r]	
				& Y\\
				PX \ar[u]^{\pi_0}\ar[d]_{\pi_1}
				& A' 
				\ar[d]|(.375)\hole|(.45)f|(.525)\hole
				\ar[l]_{h_X^P} \ar[r]^{h_Y^P}	
				& PY\ar[u]_{\pi_0}\ar[d]^{\pi_1}\\
				X
				& A
				\ar[l]^{\phi_X} \ar[r]_{\phi_Y}	
				& Y\mathrlap.
			}
		\end{aligned}
	}
	Since the $\pi_j$ are quasi-isomorphisms,
	three applications of \Cref{thm:Tor-quism}
	let us set
	\[
	\Tor_{f}(u,v;h_X,h_Y) 
	\ceq
	\Tor_{f}(\pi_1,\pi_1) \o \Tor_{\id}(\pi_0,\pi_0)\- \o \Tor_{\id}(u,v)\mathrlap.
	\]
    When \eqref{eq:Tor-DGA-functoriality-squares} 
    genuinely commutes, 
    this composite reduces to the original $\Tor_f(u,v)$.
\end{proof}

To make diagrams to come fit the margins, 
we introduce an abbreviation convention.

\begin{notation}\label{def:suppression}
	Given \DGA maps $X \from A \to Y$,
	functors $F,G,F',G'\: \Algs \lt \Algs$,
	and natural transformations
	$F \lt G$,\, $F' \lt G$,\, 
	$\phi\:F \lt F'$,
	and $\psi\:G \lt G'$
	such that the two compositions $F \lt G'$ are equal,
	we make the abbreviations
	\[
	\defm{\Tor_{FA}} \ceq \Tor_{FA}(FX,FY),
	\qquad\mnn\qquad\mnn\qquad\mnn
	\defm{\Tor_{FA}(GX)} \ceq \Tor_{FA}(GX,GY)\mathrlap,
	\]
	\vspace{-1.5em}
	\[
	\defm{\Tor_\phi} \ceq \Tor_\phi(\phi,\phi)\:
	\Tor_{FA} \lt \Tor_{F'A},
	\qquad\quad
	\defm{\Tor_\phi(\psi)} \ceq \Tor_\phi(\psi,\psi)\:
	\Tor_{FA}(GX) \lt \Tor_{F'A}(G'X)\mathrlap.
	\]
	Accordingly, in diagrams involving 
	functors applied to the triple $X \from A \to Y$,
	we will sometimes omit the ``$A$--$Y$'' squares	
	when they are symmetric with the ``$A$--$X$'' squares.
\end{notation}

\section{SHC-algebras}\label{sec:SHC}

A commutative \DGA $A$ is one for which the multiplication 
$\mu\: A \ox A \lt A$ is itself a \DGA homomorphism.
Cohomology rings are of this sort, 
and a large part 
of why homotopy theory is so much more tractable
over a field $\kk$ of characteristic $0$ is that
there are functorial \CDGA models for cochains.
For other characteristics
this is not the case~\cite[Thm.~7.1]{borel1951leray}, 
but we can weaken the requirement
by asking only that $\mu$ extend to an \Ai-algebra map.
Munkholm's product is defined in terms of such a structure,
as first considered by Stasheff and Halperin.

\begin{definition}[Stasheff--Halperin~{\cite[Def.~8]{halperinstasheff1970}}]\label{def:WHC}
	We refer to a \DGA $A$
    equipped with a \DGC map $\defm{\Phi_A}\: \B(A \ox A) \lt \B A$ 
	such that the composition $t_A \o \Phi \o \desusp_{A \ox A}\: 
	\ol{A \tensor A} \lt \ol A$
	is 
	the multiplication $\muA\: A \ox A \lt A$
	as a \defd{weakly homotopy commutative} (\textcolor{RoyalBlue}{\WHC}\defd{-})\defd{algebra}.
	Given two {\WHCA}s $A$ and $Z$,
	a \textcolor{RoyalBlue}{\WHC}\defd{-algebra map}
	from $A$ to $Z$ is 
	a \DGC map $g\: \B A \lt \B Z$
	such that there exists a \DGC homotopy between 
	the two paths around
	the square
\quation{\label{eq:SHC-def-square}
		\begin{aligned}
	\xymatrix@C=1.75em@R=3em{
		\mathllap\B(A\ox A)  
		\ar[d]_{g \,\,\mn\T\,\mnn g}		
		\ar[r]^(.565){\Phi_A}
		&
		\B A
		\ar[d]^{g}								\\
		\mathllap\B(Z \ox Z)
		\ar[r]_(.565){\Phi_Z}			
		\ar@{}[r]_{{\vphantom{g}}}				&
		\B Z\mathrlap.
	}
		\end{aligned}
}
If the square commutes on the nose,
we say the \WHCA map $g$ is \defd{strict}.
A \WHCA is called a \emph{strongly homotopy commutative} (\SHC-)algebra
if additionally $\Phi_A$
satisfies three axioms ensuring unitality, commutativity,
and associativity up to homotopy,
which for our purposes we need not spell out here.\footnote{\
		Stasheff--Halperin call
		the bare map itself, without any axioms,
		a \emph{strongly homotopy commutative algebra} 
		structure, but 
        we follow Munkholm in repurposing the term for 
        the more restrictive notion.
		}	
\end{definition}

\bex\label{thm:CGA-SHC}
If $A$ is a \CDGA, 
then the morphism $\Phi = \B\muA\: \B(A \ox A) \lt \B A$ makes $A$ 
an \SHCA.
The cohomology ring $\H(X;\kk)$ of a simplicial set is of this type,
and will always come considered with this \SHCA structure.
If $\rho\: A \lt B$ is a map of \CDGAs,
then $\B\rho$ is an \SHCA map.
\eex




\begin{theorem}[{\cite[Prop.~4.7]{munkholm1974emss}}]\label{thm:SHC-cochain}
	Let $X$ be a simplicial set and $\kk$ any ring.
	Then the normalized cochain algebra $\defm{\C(X)} = \C(X;\kk)$
	admits an \SHCA structure $\Phi_{\C(X)}$,
	and this structure is \defd{strictly natural} in the sense
	that given a map $f\: Y \lt X$ of simplicial sets,
	$\B \C(f)\: \B \C(X) \lt \B \C(Y)$
	renders the square \eqref{eq:SHC-def-square} commutative on the nose.
\end{theorem}

This natural \SHC structure on cochains 
is a reinterpretation of the classical Eilenberg--Zilber theorem;
it is only verifying the homotopy-associativity axiom that requires
substantial additional work. 
The most general known class of examples of {\SHCA}s is the following.

\bdefn\label{rmk:operads}
A \defd{homotopy Gerstenhaber algebra} (\defm{\HGA})~\cite{gerstenhabervoronov1995}
is a module~\cite{mccluresmith2003}
over the $E_2$-operad $F_2 \ms X$,
a filtrand of the \defd{surjection operad}~\cite{bergerfresse2004operad}
of interval-cut operations on cochains 
and a quotient of the second filtrand $F_2 \ms E$ of the \DG-operad $\ms E$ 
associated to the classical Barratt--Eccles simplicial operad%
~\cite{bergerfresse2004operad}.
Similarly, an \defd{extended homotopy Gerstenhaber algebra} is a module over 
a certain suboperad of the $E_3$-operad $F_3 \ms X$%
~\cite{franz2019shc}.
\edefn

\begin{theorem}[Franz~{\cite{franz2019shc}}]\label{thm:Franz-SHC}
An \HGA $A$ admits a \WHCA structure $\Phi_A$
satisfying the unitality and associativity axioms
and strictly natural in maps of \HGAs.
If $A$ is an \emph{extended} \HGA,
then $\Phi_A$ also satisfies the commutativity axiom.
\end{theorem}

\section{The products on cohomology and cochains}\label{sec:products}
To motivate Munkholm's product,
it is easiest to first 
follow him in interpreting 
the classical products on $\Tor_{\C B}(\C X, \C E)$
and $\Tor_{\H B}(\H X, \H E)$
in terms of the canonical \SHCA structures,
going into a bit more detail than he did.

The latter is the easier, so we start there.
Given \DGAs $R_0$, $R_1$ and right and left \DG $R_i$-modules 
$M_i$ and $N_i$ respectively,  
there is a classically defined external product~%
\cite[p.~206]{cartaneilenberg}
\[
\Tor_{R_0}(M_0, N_0) \ox\Tor_{R_1}(M_1, N_1)
\lt 
\Tor_{R_0 \ox R_1}(M_0 \ox M_1, N_0 \ox N_1)\mathrlap,
\]
functorial in all six variables in the sense that
given similarly defined $R'_i$, $M'_i$, $N'_i$
such that the squares \eqref{eq:Tor-DGA-homotopy-squares} commute,
then so does the square
\[
\xymatrix{
	\Tor_R(M,N) \ox \Tor_R(M,N) \ar[r]\ar[d]&\Tor_{R \ox R}(M \ox M,N \ox N)\ar[d]\\
	\Tor_{R'}(M',N') \ox \Tor_{R'}(M',N') \ar[r]&\Tor_{R' \ox R'}(M' \ox M',N' \ox N'),
}
\]
and given further $R''_i,M''_i,N''_i$, such squares glue.
If $R = R_0 = R_1$ is a \emph{commutative} \DGA,
then $\mu\: R' = R \ox R \lt R$ is a \DGA map,
and if $M = M_0 = M_1$ and $N = N_0 = N_1$ are themselves \DGAs,
then $\mu\: M' = M \ox M \lt M$
and $\mu\: N' = N \ox N \lt N$
make a diagram of the shape
\eqref{eq:Tor-DGA-functoriality-squares} commute,
so we may follow the external product with the map
\[
\Tor_\mu = \Tor_\mu(\mu,\mu)\:	\Tor_{R \ox R}(M \ox M, N \ox N)
\lt
\Tor_R(M, N)
\]
to obtain the classical product on Tor.
This particularly applies to $R = \H(B)$, \ $M = \H (X)$, \ $N = \H (E)$
for $X \from B \to E$ maps of spaces.

To define the product on Tor of cochain algebras,
we first recall its definition.
We begin by applying the external product for 
$(M \from R \to N) = \big(\C(X) \from \C(B) \to \C(E)\big)$.
We would like to apply $\Tor_\mu$ 
for $\mu$ the cup product, but cannot quite.
Recall that the cup product on $\C(B)$ can be seen as the composition
\[
\C(B) \ox \C(B) 
\os\evquism\lt 
(C_* B \ox C_* B)^*
\xtoo{a^*}
{\C(B \x B)}
\xtoo{\C(\D)} 
\C(B)
\]
where $\defm i\: \C(B) \ox \C(B) \lt (C_* B \ox C_* B)^*$ 
takes the tensor product $c \ox c'$ of cochains
to the operation $\s \ox \s' \lmt c(\s)c'(\s')$,
where
$\defm{a^*}\:(C_* B \ox C_* B)^* \lt \C(B \x B)$ is the dual 
of the classical Alexander--Whitney chain map $a$,
and where $\defm\D\: B \lt B \x B$ is the diagonal.
The hitch is that although~$\C(\D)$ and~$i$ are \DGA maps, $a^*$ is not.
However, the dual $\defm{\nabla^*}$ to the
Eilenberg--Zilber map $\nabla$
is a \DGA map natural in spaces,
meaning $\Tor_{\nabla^\medstar}$ is defined,
and the Eilenberg--Zilber theorem states that $\nabla$
is homotopy-inverse to $a$,
and hence a quasi-isomorphism,
so $\Tor_{\nabla^\medstar}$ is an isomorphism by \Cref{thm:Tor-quism}.
Thus one can follow it backward in 
the following composition 
(which uses the abbreviation convention of \Cref{def:suppression})
to obtain a product on Tor of cochain algebras:
\quation{\label{eq:Tor-C}
\Tor_{\C(B)} \ot \xtoo{\mathrm{external}} \Tor_{\C(B) \ox \C(B)} 
                \xtoo{\Tor_i} \Tor_{(C_*B \ox C_*B)^*}
                \xleftarrow[\sim]{\Tor_{\nabla^\medstar}} \Tor_{\C(B \x B)}
                \xtoo{\Tor_{C^\medstar(\D)}} \Tor_{\C(B)}
                \mathrlap.
		}
When $E \lt B$ is a Serre fibration and $\pi_1(B)$ acts trivially on the cohomology
of the homotopy fiber,
this product is sent to the cup product 
under the isomorphism with $\H(X \x_B E)$%
~\cite[Corollary 7.18]{mcclearyspectral}%
\cite[Cor.~3.5]{gugenheimmay}%
\cite[Prop.~3.4]{smith1967emss}%
\cite[Thm.~A.27]{carlsonfranzlong}.%
\footnote{\ 
	No source the author knows demonstrates this in full detail,
	but McCleary reduces it to an exercise,
	and Carlson--Franz spell out some of the steps to this exercise.
}

To realize this product as an instance of Munkholm's,
we must use his \SHCA structure on $\C$.
He notes~\cite[2.6]{munkholm1974emss}
that the Eilenberg--Zilber theorem
can be restated as the claim
\[
	\xymatrix@C=1.5em@R=1.5em{
 		\ar@(dl,ul)[]^{h} 
        & 
        \!\!\!\!\!\!\!\!\!\!\!\!
            \C(X \x X) 	
            \ar@<0.3ex>[r]^(.4){\nabla^\medstar}
            &
            {(C_* X \ox C_* X)^*}
            \ar@<0.3ex>[l]^(.6){a^\medstar} 
		}
\]
is a trivialized extension as in \Cref{def:TEX}
for a certain cochain homotopy $h$.
Then \Cref{thm:promotion} promotes $a^*$ to a twisting cochain 
$t^{a^\medstar}\: \B(C_*X\ox C_*X)^* \lt \C(X \x X)$
such that $\nabla^* \o t^{a^\tinystar} = t^{(C_*X\ox C_*X)^*}$,
with
an associated \DGC map
$g_{t^{a^\tinystar}}\: \B(C_*X\ox C_*X)^* \lt \B\C(X \x X)$.
Munkholm's \SHCA stucture $\Phi_{\C(X)}$ from \Cref{thm:SHC-cochain}
is the composite
\[
\B(\C X \ox \C X) \xtoo{\B i}
\B(C_*X\ox C_*X)^* \xtoo{g_{t^{a^\tinystar}}} 
\B\C(X \x X) \xtoo{\B \C(\D)}
\B\C(X)\mathrlap.
\]

Agreeing to write 
$A = (C_*X \ox C_*X)^*$ for brevity,
note that the recursive prescription 
$t^{a^\tinystar} = h(t^{a^\tinystar} \cup t^{a^\tinystar}) + a^* t^{A}$
of \Cref{thm:promotion}
means that $t^{a^\medstar}|_{\B_1 A}$ 
is just $a^* t^{A}|_{\B_1 A}$,
so that if we write $\desusp\: \ol A \simto \B_1 A$
for the inverse to $t^{A}|_{\B_1 A}$,
then 
\[
	a^* 
	=
    t^{a^\medstar} \o \desusp
    =
    t^{\wt A} \o g_{t^{a^\tinystar}} \o \desusp
    =
    \e \o t_{\B \wt A} \o g_{t^{a^\tinystar}} \o \desusp 
    =
    \e \o \W g_{t^{a^\tinystar}} \o t_{\B A} \o \desusp
    =
    \e \o \W g_{t^{a^\tinystar}} \o i_A\mathrlap,
\]
by the definitions of the tautological twisting cochains,
where $i_A$ is the \DG module section of $\e$ from \Cref{def:TEX-WBA}.
Thus, in cohomology,  
$\H(\nabla^*)\- = \H(a^*) = \H(\e) \o \H(\W g_{t^{a^\tinystar}}) \o \H(\e)\-$.
%

%

The following squares then allow us to substitute the last 
three maps in \eqref{eq:Tor-C}, along the top, with the composition $\Tor_{\W \Phi}$
along the bottom:
%
\[
\xymatrix@R=4.5em@C=3.5em{
	\Tor_{\C (B)\ox \C (B)} 			\ar[r]^{\Tor_\evquism}
	&\Tor_{(C_* B \ox C_* B)^*}
	&\Tor_{\C(B \x B)} 			\ar[r]^(.5375){\Tor_{\C(\D)}} 
	\ar[l]_(.475){\Tor_{\nabla^*}}^\sim
	&\Tor_{\C(B)}
	\\
	\Tor_{\W\B(\C B\ox \C B )} 	\ar[r]_{\Tor_{\W\B\evquism}}
	\ar[u]^{\Tor_\e}_\vertsim
	&\Tor_{\W\B(C_* B \ox C_* B)^*}	
	\ar[u]^{\Tor_\e}_\vertsim
	&\Tor_{\W\B\C(B \x B)}			\ar[r]_(.5375){\Tor_{\W\B\C(\D)}}
	\ar@{<-}[l]^(.475){\Tor_{\W g_{\mn t^{a^{\mnn\smash\tinystar}}}}}
	\ar[u]_{\Tor_\e}^\vertsim
	&\Tor_{\W\B\C(B)}\mathrlap.		\ar[u]_{\Tor_\e}^\vertsim
}
\]

%
%
The preceding diagram is compressed to the last square in the following:
\[
	\begin{aligned}
		\xymatrix@R=4.5em{
	        \Tor_{\C(B)} \ox \Tor_{\C(B)}
				\ar[r]^(.525){\mr{external}}
		&\Tor_{\C(B) \ox \C(B)} 
				\ar@{=}[r]
		&\Tor_{\C(B) \ox \C(B)} 
				\ar[r]
		&\Tor_{\C(B)}
		\\
	        \Tor_{\W\B\C(B)}\ox\Tor_{\W\B\C(B)}
				\ar[r]_(.5375){\mr{ext.}}
                                \ar[u]^{\Tor_\e \ox \Tor_\e}_\vertsim
		&\Tor_{\W\B\C(B) \ox \W\B\C(B)} 
		\ar@{<-}[r]_(.5625){\Tor_\psi}^(.5625)\sim
				\ar[u]^{\Tor_{\e \ox \e}}_\vertsim
		&\Tor_{\W\B(\C B  \ox \C B )} 
				\ar[r]_(.575){\Tor_{\W\Phi}}
				\ar[u]_{\Tor_\e}^\vertsim 
		&\Tor_{\W\B\C(B)}\!\mathrlap,
				\ar[u]_{\Tor_\e}^\vertsim
		}
	\end{aligned}
\]
which commutes by naturality of the external product and by \Cref{def:psi},
so finally we see the classical product on $\Tor_{\C(B)}$ 
realized as the composite along the other three outer edges 
of the big rectangle.

This same construction evidently applies to $\H(X) \from \H(B) \to \H(E)$ 
with $\Phi_{H^\medstar(B)} = \B\mu_{H^\medstar(B)}$ and so on.
More generally,
this construction applies for a triple $X \from A \to Y$
of \DGA maps such that the induced maps $\BX \from \BA \to \BY$
are strict \WHCA maps; in this more general case there is not necessarily any
composite along the top anymore, so the product
can only be defined as the composite along the other edges,
$\Tor_\e \o \Tor_{\W\Phi} \o \Tor_\psi\- \o \,\mr{ext.}\mnn \o (\Tor_\e\ot)\-$.

\section{Munkholm's product}\label{sec:product}

Munkholm now generalizes the product on Tor
from the case of \DGA maps $X \from A \to Y$
inducing strict \WHCA maps
to the case of a general triple $\B X \from \B A \to \B Y$ of \SHCA maps,
which are not induced by \DGA maps. 
Thus we are assuming the following homotopy-commutative squares of \DGC maps.

\begin{equation}\label{eq:SHC-map}
\begin{aligned}
\xymatrix@C=1.75em@R=4em{
	\B(X \ox X)\ar[d]_{\Phi_X}& 
	\B(A \ox A)\ar[d]|(.45)\hole|{\Phi_A}|(.55)\hole
	\ar[l]_{\xi \,\T\, \xi}
	\ar[r]^{\upsilon \,\T\, \upsilon}	
	& \B(Y \ox Y)\ar[d]^{\Phi_Y} \\
	\B X	& \B A			\ar[l]^{\defm\xi}\ar[r]_{\defm\upsilon}	
	& \B Y
}
\end{aligned}
\end{equation}
Taking $\W$ of the diagram, one has a homotopy-commutative diagram of \DGAs,
inducing a map $\Tor_{\W\B(A \ox A)} \lt \Tor_{\W\BA}$ 
by the technique of \Cref{thm:homotopy-Tor-map}.
For the special cases of cochain algebras and cohomology rings,
this map reduces to $\Tor_{\W\Phi}$.

Working backward emulating the last diagram for the product on $\Tor_{\C(B)}$,
one wants to pass from $\Tor_{\W\BA \ox \W\BA}$
to $\Tor_{\W\B(A \ox A)}$
using $\Tor_\psi = \Tor_\psi(\psi,\psi)$, 
but one cannot do so in general
as $\psi\:\W\B(- \ox -) \lt \W\B(-) \ox \W\B(-)$
is natural only in pairs of \DGA maps, 
and a \DGC map like $\W(\xi \T \xi)$
cannot be assumed to be of the form $\W\B f$ 
for a \DGA map $f\: A \ox A \lt X\ox X$.
To work around this, Munkholm employs 
the following noncommutative diagram,
in which we are omitting $A$--$Y$ squares:
\quation{\label{eq:product-setup}
	\begin{aligned}
\xymatrix@C=1.5em@R=2.5em{
	\W \B (A \ox A)		
	\ar[d]_{\W(\xi \,\T\, \xi)}
	\ar^(.45)\psi[r]
	&
	\W \B A \ox \W \B A	
	\ar[d]|(.425)\hole|(.525)\hole
	|{\W\xi \ox \W\xi}
	\ar@{=}[r]	
	&
	\W \B A \ox \W \B A	
	\ar[d]^{(\e\, \W \xi)^{\ox 2}}
	\\
	\W \B(X \ox X)		\ar_(.45)\psi[r]
	&
	\W \B X \ox \W \B X	\ar[r]_(.6){\e \ox \e} 
	& 
	X \ox X
    \mathrlap.
}
	\end{aligned}
}
Although the left square does not commute,
the right square does by definition,
and the large outside rectangle 
does too,
for $\e\ot \o (\W\xi)\ot \o \psi = \e \o \W(\xi \T \xi)$ by \Cref{thm:epsilon-T}
and $\e\ot \o \psi = \e$ by \Cref{def:psi}.
Thus we may use the right rectangle
to apply $\Tor_{\id}(\e \ot )$,
and then use the big rectangle
to apply the inverse of
$\Tor_{\psi}(\e)$:
\quation{\label{eq:substitute}
	\us{(\W \B A)^{\ox 2}}\Tor\big(\mnn(\W \B X)^{ \ox 2}\big)
	\xtoo[\sim]{\Tor_{\id}(\e\ot)}
	\us{(\W \B A)^{\ox 2}}\Tor(X\ot)
	\xleftarrow[\sim]{\Tor_\psi(\e)}
	\us{\W \B (A\ot)}\Tor\big(\W\B(X\ot)\mnn\big)
	\mathrlap.
}
In case the given \DGC maps are induced
by \DGA maps $X \from A \to Y$,
the left does square commute by naturality, so 
$\Tor_\psi = \Tor_\psi(\psi)$ makes sense and 
$\Tor_\psi(\e) = \Tor_\psi(\e\ot\o\psi) = \Tor_{\id}(\e\ot)\o\Tor_\psi$,
and hence \eqref{eq:substitute} reduces to the $\Tor_\psi\-$
appearing in the construction of the previous section.

All told, one finally gets the following composite:

	\quation{\label{eq:full-product}
		\Big(\,\us{\W\B A}\Tor\,\Big){}\ot
		\xtoo{\!\mr{ext.}\!\!}
     \!\!\!\!
		\us{(\W\B A)\ot}{\Tor}
     \!\!\!\!
		\xtoo{\!\!\!\Tr\id(\e\ot)\!\!\!\!\!\!}
     \!\!\!\!
     \us{\W(\B A)\ot\vphantom{x^h}}{\Tor}\!\!\!\!\!(X\ot)
		\xleftarrow{\!\!\!\!\!\Tr{\psi}(\e)\!\!\!\!}
	\!\!\!\! \! 
		\us{\W\B(A\ot)}{\Tor}
	\!\!\!\!
		\xtoo{\!\!\!\Tr{\id}(\W\Phi)\!\!\!\!\!\!}
	\!\!\!\!
		\us{\W\B(A\ot\vphantom{x^h})}\Tor\!\!\!\!\!(\W\B X)
		\xtoo{\!\!\!\Tr\id(\pi_0)\!\!\!}
	\!\!\!\!\!
		\us{\W\B(A\ot\vphantom{x^h})}\Tor\!\!\!\!\!(P\W\B X) 
		\xleftarrow{\!\!\!\Tr{\W\Phi}(\pi_1)\!\!\!}
		\us{\W\B A}{\Tor}\mathrlap.
	}

\section{The reformulation}\label{sec:CGA}

It would be helpful to have a description of the substitute
$\Tor_{\psi}(\e)
\- \o \Tor_{\id}(\e\ot)$
of \eqref{eq:substitute}
that behaves uniformly in the three variables of Tor.
We accomplish this by
replacing it with $\Tor_{\W\mnn\nabla} \o \Tor_\g \-$.

\begin{theorem}\label{thm:WDg}
	Let $A_1$ and $A_2$ be \DGAs.
	Then the composition 
	\[	
		\W(\B A_1 \ox \B A_2) \xtoo{\W\mnn\nabla}
		\W\B(A_1 \ox A_2) \xtoo{\psi}
		\W\B A_1 \ox \W\B A_2
	\]
	agrees with $\g$ from \Cref{def:shuffle}.
\end{theorem}

\begin{proof}
We will show 
the twisting cochains
$\B A_1 \ox \B A_2 \lt \W\B A_1 \ox \W\B A_2$ 
associated to $\g$ and $\psi \o \W\mnn\nabla$ are equal.
The former, $\gamma t_{\B A_1 \ox \B A_2}$,
is $t_{\B A_1} \ox \h\e + \h\e \ox t_{\B A_2}$
by \Cref{def:shuffle},
whereas by naturality of $t_{({-})}$ the latter is 
\[
	\psi \o \W\mnn\nabla \o t_{\B A_1 \ox \B A_2} = 
	 \psi \o t_{\B(A_1 \ox A_2)}\o \nabla =
	 t_\psi \o \nabla
\mathrlap.
\]
Because $t_\psi$ is given
via the recursive prescription
$t_\psi = h_\psi(t_\psi \cup t_\psi) +  (i_1 \ox i_2) t^{A_1 \ox A_2}$
of \Cref{thm:promotion},
the restriction of the twisting cochain $t_\psi \nabla$ to
$\nabla\-\big( \B_1 (A_1 \ox A_2)\big)$ 
is given by
\[		
				(i_1 \ox i_2)
				(t^{A_1} \ox \h\e + \h\e \ox t^{A_2})
			 = t_{\B A_1} \ox \h\e 
				+ 
				\h\e \ox t_{\B A_2}\mathrlap,
\]
agreeing with $\g t_{\B A_1 \ox \B A_2}$.
We must check that they agree on all of 
$\B A_1 \ox \B A_2$.

We begin with $\B A_1 \ox \kk$,
on which $\nabla$ restricts to an isomorphism
$\B A_1 \ox \kk \lt \B(A_1 \ox \kk)$.
The image of $\B_1(A_1 \ox \kk)$
under $(i_1 \ox i_2)t^{A_1 \ox A_2}$
lies in $\W\B A_1 \ox \kk$,
so one can write 
\[
	(i_1 \ox i_2)t^{A_1 \ox A_2}\nabla =
    i_1t^{A_1} \ox \h
	\: \B A_1 \ox \kk \lt \W\B A_1 \ox \kk
	\mathrlap.
\]
The homotopy $h_j$ from \Cref{def:TEX-WBA} vanishes on $\kk = \W_0 \B A_j$
for $j \in \{1,2\}$,
so one has $h_\psi = h_1 \ox {\id} + i\e \ox h_2 = h_1 \ox \id$
on $\W\B A_1 \ox \kk$,
one sees from the recursive prescription 
$t_\psi|_{B_{\geq 2}(A_1 \ox \kk)} = h_\psi(t_\psi \cup t_\psi)$
that this process is effectively the same as that defining 
the cochain $t^{i_{A_1}} = t_{\B A_1}$ 
of \Cref{def:TEX-WBA}, 
but with added inert ``$1$'' tensor factors.
Thus $t_\psi$ agrees with
$t_{\B A} \ox \h_{\W \B A}$ on $\B A \ox \kk$.
The proof for $\kk \ox \B A_2$ is symmetric,
noting that $i_{2}\e_{A_2}(1) = 1 \in \W_0\B A_2$.

It remains to see $t_\psi \nabla$ vanishes on 
$\B_{\geq 1} A_1 \ox \B_{\geq 1} A_2$.
Start small, letting
$b_1 \ox b_2 \in \B_1 A_1 \ox \B_2 A_2$
be given;
we are to evaluate
\[
	t_\psi \nabla (b_1 \ox b_2) =
	(h_1 \ox {\id} + i_1\e \ox h_2)\mu_{\W\B A_1 \ox \W \B A_2}
	(t_\psi \ox t_\psi)
	\D\nabla(b_1 \ox b_2)\mathrlap.
\]
Note that
$\D_{\B(A_1 \ox A_2)}\nabla = (\nabla \ox \nabla)\D_{\B A_1 \ox \B A_2}$.
By \Cref{lem:symmetry},
we know $(\nabla \ox \nabla)\D(b_1 \ox b_2)$
is the sum of 
$\nabla(b_1 \ox 1) \ox \nabla(1 \ox b_2)$
and $(-1)^{|b_1||b_2|}\nabla(1 \ox b_2) \ox \nabla(b_1 \ox 1)$.
Now $t_\psi$ is defined to be
$(i_1\ox i_2)t^{A_1 \ox A_2}$
on 
$\B_1(A_1 \ox \kk)$
and $\B_1(\kk \ox A_2)$,
which respectively contain $\nabla(b_1 \ox 1)$
and $\nabla(1 \ox b_2)$,
and the image of 
$i_j\: A_j \lt \W\B A_j$
lies in $\W_1 \B_1 A_j \leq S_0 \leq \ker h_j$
by \Cref{def:TEX-WBA}, so 
(suppressing ``1'' and ``2'' subscripts out of space considerations)
\eqn{
	(h \ox {\id} + i\e \ox h)\mu
	\big(\mn(i \ox i)t^{A_1 \ox A_2}
	\big){}\ot
	\big(\mnn\nabla(b_1 \ox 1) \ox \nabla(1 \ox b_2)\mn\big)
	&=
		(h \ox {\id} + i\e \ox h)\mu
		\big(\mn(ib_1 \ox i1) \otimes (i1 \ox ib_2)\mn\big)
	\\&=
	(h \ox {\id} + i\e \ox h)(ib_1 \ox ib_2)
	\\&=
	hib_1 \ox ib_2 + i\e i b_1 \ox hib_2 = 0\mathrlap,
	}
and similarly for the other summand.

%

Now let $b_1 \ox b_2$ lie in $\B_m A_1 \ox \B_n A_2$
and suppose inductively that we know $t_\psi$
vanishes on $\B_p A_1 \ox \B_q A_2$
for pairs $(p,q) \neq (m,n)$ with 
$1 \leq p,q$ and $p \leq m$ and $q \leq n$.
We know $t_\psi\nabla 
= 
h_\psi \mu
(t_\psi\nabla \ox t_\psi\nabla)\D
$,
where $t_\psi$ vanishes on $\B_0(A_1 \ox A_2)$
and by \Cref{thm:BA-structure}, 
$\D$ sends $\B_m A_1 \ox \B_n A_2$
to the sum of terms in 
$(\B_p A_1 \ox \B_q A_2) \otimes (\B_{p'} A_1 \ox \B_{q'}A_2)$
with $p+p' = m$ and $q+q' = n$.
By the induction hypothesis,
the only terms of $\D(b_1 \ox b_2)$
not necessarily annihilated by 
$h_\psi \mu
(t_\psi\nabla \ox t_\psi\nabla)$
are those lying in $(\B_m A_1 \ox \kk) \otimes (\kk \ox \B_n A_2)$
and $(\kk \ox \B_n A_2) \otimes(\B_m A_1 \ox \kk)$,
to wit,
 $b_1 \ox 1 \ox 1 \ox b_2$ and 
 $(-1)^{|b_1||b_2|} 1 \ox  b_2 \ox b_1 \ox 1$.
But 
by \Cref{lem:symmetry} and the fact $t_\psi$ is of degree $1$,
we have
\eqn{
	\mu(t_\psi \ox t_\psi)\D\nabla(b_1 \ox b_2) 
&=
	\mu(t_\psi \ox t_\psi)
	\big(\mnn
		\nabla(b_1 \ox 1) \ox \nabla(1 \ox b_2) 
		+ (-1)^{|b_1||b_2|}
		\nabla(1 \ox b_2) \ox \nabla(b_1 \ox 1)
	\mnn\big)
\\ 
&=
	\mu\Big(\mn
		(-1)^{|b_1|}(t_{\B A_1}b_1 \ox 1) \otimes 
			(1 \ox t_{\B A_2}b_2) 
\\
&
	\phantom{
				=(-1)^{|b_1|}(t_{\B A_1}b_1 \ox 1) 
				\otimes b \otimes b \,
			}
		+ (-1)^{(|b_1|+1)|b_2|} 
		(1 \ox t_{\B A_2} b_2) \otimes 
			(t_{\B A_1}b_1 \ox 1)
	\mnn\Big)
\\ 
&=
	\big(\mn
		(-1)^{|b_1|} + 
		(-1)^{(|b_1|+1)|b_2|\ +\ (|b_1| +1)(|b_2|+1)}
	\big)\,\mn
	\mu
	\big(\mn
		(t_{\B A_1}b_1 \ox 1) \otimes
		(1 \ox t_{\B A_2}b_2)
	\mn\big) 
\\
&= 0
\mathrlap.\qedhere
}
%
%
%
\end{proof}

Thus we can replace the second two maps 
in \eqref{eq:full-product}
to obtain the more tractable product whose properties are explored in
the predecessor~\cite{carlson2022munkholm}.

\begin{corollary}\label{thm:Jeff-equiv-Munkholm}
Given \WHCA maps and homotopies as
in \eqref{eq:SHC-map},
the product \eqref{eq:full-product-v2}
can be equivalently expressed as the composite	
	\quation{\label{eq:full-product-v2}
		\Big(\,\us{\W\B A}\Tor\,\Big){}\ot
		\xtoo{\!\mr{ext.}\!}
	\!\!\!
		\us{(\W\B A)\ot}{\Tor}
	\!
		\xleftarrow[\sim]{\!\Tr\g}
	\!
		\us{\W(\B A)\ot}{\Tor}
	\!
		\xtoo{\!\mn\Tr{\W\nabla}\!\!}
	\!
		\us{\W\B(A\ot)}{\Tor}
	\!\!\!
		\xtoo{\!\!\!\Tr{\id}(\W\Phi)\!\!\!\!}
	\!\!\!
	\us{\W\B(A\ot\vphantom{x^h})}\Tor\!\!\!\!\!(\W\B X)
		\xtoo{\!\!\!\Tr{\id}(\pi_0)\!\!\!\!}
	\!\!\!
	\us{\W\B(A\ot\vphantom{x^h})}\Tor\!\!\!\!\!(P\W\B X)
		\xleftarrow[\sim]{\!\!\mn\Tr{\W\Phi}(\pi_1)\!\!\!}
		\us{\W\B A}{\Tor}\mathrlap.
	}
\end{corollary}

\begin{proof}
First note that
these maps of Tors are well-defined by naturality of $\g$ and $\W\mnn\nabla$,
and $\Tor_\g$ is invertible by \Cref{thm:Tor-quism}
since $\g$ is a quasi-isomorphism by \Cref{def:shuffle}.
Now recall from \Cref{thm:WDg} that $\psi \o \W\mnn\nabla = \g$
and from \Cref{def:psi} that $\e\ot \ox \psi = \e$,
so that 
\[
\e\ot \o \g = \e\ot \o \psi \o \W\mnn\nabla = \e \o \W\mnn\nabla\mathrlap.
\]
Hence the following diagram commutes,
and symmetrically for the $A$-$Y$ squares:
\[
\xymatrix@C=1.5em@R=3.75em{
	\W(\B A)\ot 
	\ar@/^1.25pc/[rrr]^\g
	\ar@{}[r]_(.475){\W\mnn\nabla}
	\ar@{}[r]^(.125){}="a"^(.3){}="b" 
	\ar "b";[r]	
	\ar[d]					&
	\W\B(A\ot)
	\ar[r]_(.5){\psi}		
	\ar[d]					&
	(\W\B A)\ot
	\ar@{=}[r]				
	\ar[d]					&
	(\W\B A)\ot	
	\ar[d]				
	\\
	\W(\B X)\ot
	\ar@{}[r]^(.475){\W\mnn\nabla}
	\ar@{}[r]^(.2){}="a"^(.3){}="b" 
	\ar "b";[r]	
	\ar@/_1.25pc/@{{}{ }{}}[rrr]^(.05){}="c"
	\ar@/_1.25pc/"c";[rrr]_\g&
	\W\B(X\ot)
	\ar[r]^(.525)\e				&
	X\ot						&
	(\W\B X)\ot
	\ar[l]_(.5375){\e\ot}
	\mathrlap.
}
\]
%
Thus we find
$
\Tor_{\id}(\e\ot) \o \Tor_\g = \Tor_\psi(\e) \o \Tor_{\W\mnn\nabla}
$.
Replacing $\Tor_{\psi}(\e)\- \o \Tor_{\id}(\e\ot)$ 
in \eqref{eq:full-product} with $ \Tor_{\W\mnn\nabla} \o \Tor_\g \-$ 
yields \eqref{eq:full-product-v2}.
\end{proof}

\counterwithin{figure}{section}
\counterwithin{theorem}{section}
\counterwithin{lemma}{section}
\counterwithin{corollary}{section}
\counterwithin{proposition}{section}
\counterwithin{equation}{section}
\counterwithin{notation}{section}
\numberwithin{remark}{section}

\section{The product on the two-sided bar construction}\label{sec:comparison}
In this section we show our product on Tor is induced
by the product on the two-sided bar construction
due to Franz~\cite[App.~A]{carlsonfranzlong},
which we have used in previous work.

\begin{definition}
	Given a \DGC $C$, \DGAs $X$ and $Y$,
	and twisting cochains $\tau^X\: C\lt X$ and $\tau^Y\: C \lt Y$,
	the \defd{twisted tensor product}
	$\defm{X \ox_{\tau^X} C \ox_{\tau^Y} Y}$
	 is the complex with underlying graded $\kk$-module $X \ox C \ox Y$
	equipped with the differential given as the the sum of the tensor differential
	and the two operations 
\eqn{
		(\mu_X \ox \id_C)(\id_X \ox \tau^X \ox \id_{\BA})(\id_X \ox \D_{\BA}) \ox \id_Y
	\:&
		x \ox c \ox y \lmt \pm x \. \tau^X(c_{(1)}) \ox c_{(2)} \ox y
	\mathrlap,\\
		-\id_X \ox {}(\id_{\BA} \ox \mu_Y)(\id_{\BA} \ox \tau^Y\upsilon \ox \id_Y)(\D_{\BA} \ox \id_Y)
	\:&
		x \ox c \ox y \lmt \pm x \ox c_{(1)} \ox \tau^Y(c_{(2)})\. y
	.\footnotemark
}
    \footnotetext{\ 
        The signs $\pm$ on the right-hand sides 
        are determined by the definition of the map 
        on the left by the Koszul rule.
    }
Given a span $\B X \os{\defm\xi}\from \BA \os{\defm\upsilon}\to \B Y$ of \DGC maps,
the \defd{two-sided bar construction}
is the twisted
tensor product $\smash{\dsp\defm{\B(X,A,Y)} \ceq X \ox_{t^X \xi} \B A \ox_{t^Y \upsilon} Y}$.
\end{definition}

Twisted tensor products exhibit functoriality with respect to
commutative diagrams in which $l_X$, $l_Y$ are \DGA maps and $g$ a \DGC map, thus%
~\cite[Lem.~1.20]{carlsonfranzlong}:
\quation{%
\label{eq:twisted-tensor-map}
\begin{aligned}
\xymatrix{
\ar[d]_{l_X}
	X' & \ar[l] C' \ar[r] \ar[d]|(.475)\hole|g|(.525)\hole & Y' \ar[d]^{l_Y}\\
	X & \ar[l] C \ar[r] & Y\mathrlap.
	}
\end{aligned}
}
See Carlson--Franz~\cite[Def.~1.16 \emph{et seq.}]{carlsonfranzlong} 
for much more detail and 
Carlson~\cite{carlson2022collapses} for the history of these notions.

\begin{definition}[{Wolf~{\cite[p.~322]{wolf1977homogeneous}}}]
We write $\defm{\TOR_A(X,Y)}$ for the bigraded cohomology
$\kk$-module $\H{\B(X,A,Y)}$.
\end{definition}

\nd This is reasonable because in case the \DGC maps are 
induced from \DGA maps as $\xi = \B x$ and $\upsilon = \B y$,
the one-sided bar construction $\B(X,A,A)$
is a proper projective $A$-module resolution of $X$
under reasonable flatness hypotheses,\footnote{\
	For example,
	it is enough that $A$ and $X$ 
	be flat over the principal ideal domain $\kk$~{\cite[after Prop.~10.19]{barthelmayriehl2014}}.
	}
and then $\B(X,A,Y)$ computes $\Tor_A(X,Y)$.
Any more specific hypotheses guaranteeing this would
complicate the statements below while needlessly excluding some cases,
so we instead directly stipulate that the bar construction
compute Tor.

\begin{definition}\label{def:bariso}
Suppose given a span
$
X \os {\defm x\vphantom{y}}\from A \os {\defm y}\to Y
$
of \DGA maps.
As $\B(X,A,A) = X \ox_{xt^A} \BA$ is exact,
given a proper, projective \DG $A$-module resolution $\defm{\Pb}$ of $X$,
there exists a map 
$P^\bul \lt \B(X,A,A) = X \ox_{xt^A} \BA$
of complexes of \DG $A$-modules over $X$,
unique up to cochain homotopy.
Applying $- \ox_A Y$ yields
a map $ P^\bul \ox_{ A}  Y \lt \B(X,A,Y)$,
inducing a map $\defm\bariso\:\Tor_A(X,Y) \lt \TOR_A(X,Y)$ in cohomology,
natural in maps $(l_X,l_A,l_Y)\: (X' \from A' \to Y') \lt (X  \from A  \to Y)$ of spans of \DGA maps
as in \eqref{eq:Tor-DGA-functoriality-squares}.
Explicitly, the following square commutes.
\[
\xymatrix@C=1em{
\Tor_{A'}(X',Y')\ar[r]^(.475)\bariso\ar[d]_{\Tor_{l_A}(l_X,l_Y)}&
\TOR_{A'}(X',Y')\ar[d]^{\H(l_X \ox \B l_A \ox l_Y)}\\
\Tor_{A}(X,Y)\ar[r]_(.475)\bariso&
\TOR_{A}(X,Y)
}
\]
In particular, a span
$\B X \os{\defm\xi}\from \BA \os{\defm\upsilon}\to \B Y$ 
of \DGC maps
induces a span
$\W\B X \os{\W\xi}\from \W\BA \os{\W\upsilon}\to \W\BY$
of \DGA maps
and we can apply this definition
to obtain 
a map 
$\defm\WBbariso\:\Tor_{\W\B A}(\W\B X,\W\B Y) \lt \TOR_{\W\B A}(\W\B X,\W\B Y)$\footnote{\ 
Munkholm simply \emph {defines} ${\Tor_A(X,Y)}$
to be $\Tor_{\W\B A}(\W\B X,\W\B Y)$ when there are no \DGA maps
$X \from A \to Y$.
}
natural in commuting triples of \DGC maps
\begin{equation}\label{eq:Gamma}
\begin{aligned}
\xymatrix@C=1.25em{
	\ar[d]_{\l_X}
	\BX' & \ar[l]_{\xi'} \BA' \ar[r]^{\upsilon'} \ar[d]|(.475)\hole|{\l_A} & \BY' \ar[d]^{\l_Y}\\
	\BX & \ar[l]^{\xi} \BA \ar[r]_\upsilon & \BY\mathrlap.
}
\end{aligned}
\end{equation}
We say the span of \DGC maps \defd{satisfies sufficient flatness conditions}
when $\WBbariso$ is an isomorphism.
\end{definition}

\brmk\label{rmk:qi}
This is \emph{a priori} a bit weaker than asking the bar construction
itself be a proper projective resolution,
but is enough to make Wolf's $\TOR_A(X,Y)$
agree with 
$\Tor_{\W\B A}(\W\B X,\W\B Y)$
for a triple of \DGC maps $\B X \from \BA \to \B Y$.
If one does have \DGA maps $ X \from A \to  Y$,
then $\Tor_\e$ is a quasi-isomorphism $\Tor_{\W\BA} \lt \Tor_A$
and $\e \ox \B\e\ox\e$ is a quasi-isomorphism 
$\B(\W\BX,\W\BA,\W\BY) \lt \B(X,A,Y)$,
so it follows $\bariso$ is a quasi-isomorphism.
\ermk

Maps of two-sided bar constructions are induced not only
from triples $(l_X,g,l_Y)$ of the form \eqref{eq:twisted-tensor-map},
but also triples $(\l_X,\l_A,\l_Y)$ of \DGC maps of the form
\eqref{eq:Gamma}
(although we claim no functoriality for such maps on the cochain level)%
~\cite[Prop.~1.26]{carlsonfranzlong}.\footnote{\ 
	They can actually be defined more generally still, using the pattern for 
	the map $\Theta$ from the penultimate section of that paper,
	which resembles the maps of \eqref{eq:Tor-DGA-homotopy-squares}.
	}
\bthm[Wolf~{\cite[Thm.~7]{wolf1977homogeneous}}]\label{thm:Gamma}
A strictly commuting diagram \eqref{eq:Gamma} of \DGC maps
gives rise to a cochain
map
\[
\defm{\B(\l_X,\l_A,\l_Y)} \: 
\B(X',A',Y')
\lt
\B(X,A,Y)\mathrlap.
\]
\nd If we have $\l_X = \B l_X$ and $\l_Y = \B l_Y$ 
for \DGA maps $l_X$ and $l_Y$,
we have $\B(\l_X,\l_A,\l_Y) = l_X \ox \l_A \ox l_Y$.
\ethm
In particular,
given {\WHCA}s algebras $A$, $X$, $Y$ 
along with \DGA maps $X \os{\defm x}\from A \os{\defm y}\to Y$ making 
\begin{equation}\label{eq:bar-prod}
\begin{aligned}
\xymatrix@C=2.5em{
	\ar[d]_{\Phi_X}
	\B(X\ox X) & \ar[l]_{\B(x \ox x)} \B(A \ox A) \ar[r]^{\B(y \ox y)} \ar[d]|(.475)\hole|{\Phi_A} & \B(Y \ox Y) \ar[d]^{\Phi_Y}\\
	\BX & \ar[l]^{\B x} \BA \ar[r]_{\B y} & \BY
}
\end{aligned}
\end{equation}
commute strictly,
there is an induced map 
$\B(\Phi_X,\Phi_A,\Phi_Y)\: \B(X\ot,A\ot,Y\ot)\lt \B(X,A,Y)$,
used by Franz~\cite[Thm.~A.1]{carlsonfranzlong} 
to define a weak product on the two-sided bar construction:
\quation{\label{eq:two-sided-bar-product}
\resizebox{6.5in}{!}{
\xymatrix@C=5em{
	\B(X,A,Y)\ot 
	\ar[r]_(.345){(2\,3\,5\,4)}^(.345){\sim}&
	\dsp
	X\ot \ox_{t^{X \ot} (\B x)\ot}  
		(\BA)\ot  \ox_{t^{Y \ot} (\B y) \ot}  Y\ot 
	\ar[r]_(.6){\id \ox \nabla \ox \id}
	&
	\dsp
	\B(X\ot,A\ot,Y\ot)
	\ar[r]_(.56){\B(\Phi_X,\Phi_A,\Phi_Y)}
	&
	\B(X,A,Y)\mathrlap,
}
}
}

\nd where the first map merely permutes tensor-coordinates:
\[
(2\ 3\ 5\ 4)\:
(x_1\ox b_1 \ox y_1) \otimes (x_2 \ox b_2 \ox y_2)
	\lmt
	(-1)^{|x_2||b_1| + |x_2||y_1| + |y_1||b_2|}
(	x_1\ox x_2) \otimes (b_1 \ox b_2) \otimes (y_1 \ox y_2)
	\mathrlap.
	\]
Now, as promised, we show that this product 
induces Munkholm's product on Tor,
assuming the two-sided bar construction computes Tor in the first place.

\begin{theorem}\label{thm:Jeff-equiv-Franz}
Assuming sufficient flatness hypotheses 
on the \DGC maps
induced from the spans 
\[
X \os {x\vphantom{y}}\from A \os y\to Y
\qquad
\mbox{and}
\qquad
X \ox X \os {\,x\ox x\vphantom{y}}\longfrom A \ox A\os {y\ox y}\longto Y \ox Y
\]
of \DGA maps,
the product \eqref{eq:two-sided-bar-product}
on the two-sided bar construction induces the product of \Cref{sec:product}
on Tor.
\end{theorem}
\begin{proof}
In \eqref{eq:comparison-squares} below,
the composite along the top row is our product on Tor
and that along the bottom is induced by Franz's product
on the two-sided bar construction.

\begin{equation}\label{eq:comparison-squares}
\begin{aligned}
\xymatrix@C=2.35em@R=3em{
  \Tor_{\W\BA}\ot \ar[d]^{\Tor_\e\ox\Tor_\e}_\vertsim\ar[r]^(.425){\EXT}
& \Tor_{(\W\BA)\ot} \ar[dr]_{\Tor_{\e\ot}}^\sim
& \Tor_{\W(\BA)\ot}\ar[l]_{\Tor_\g}^\sim\ar[r]^{\Tor_{\W\nabla}}_\sim
& \Tor_{\W\B(A\ot)}\ar[d]^\WBbariso_\vertsim\ar[dl]^{\Tor_{\e}}_\sim\ar[rrr]^{\Tor_{\W\Phi}}
&
&
& \Tor_{\W\BA}\ar[d]^\WBbariso_\vertsim	
\\ \Tor_{A}\ot	\ar[rr]_{\mr{external}} \ar[d]^{\bariso\ox\bariso}_\vertsim
&
& \Tor_{A\ot} \ar[dr]_\sim^{\bariso}		
& \TOR_{\W\B(A\ot)}  \ar[d]^{\e\ox\B\e\ox\e}
\ar[rrr]^{\H(\W\Phi \ox \B\W\Phi \ox \W\Phi)}_{\H\B(\B\W\Phi,\B\W\Phi,\B\W\Phi) }
&
&
& \TOR_{\W\BA}  \ar[d]^{\mathrlap{\e\ox\B\e\ox\e}}		
\\ \TOR_A\ot \ar[rrr]_{\H((\id \ox \nabla \ox \id)(2\,3\,5\,4))}
&
&
& \TOR_{A\ot}\ar[rrr]_(.515){\H\B(\Phi,\Phi,\Phi)} 
&
&
& \TOR_A
}
\end{aligned}
\end{equation}

\nd We will have proven the theorem if we can show \eqref{eq:comparison-squares}
is commutative.


\nd 1. The upper-left trapezoid commutes by naturality of the external product.

\nd 2. The upper triangle, in truth a diamond, can be subdivided as 
\[
\xymatrix@C=-1em{
	& \Tor_{\W(\BA)\ot} \ar[dr]^{\Tor_{\W\nabla}}\ar[dl]_{\Tor_\g}\\
\Tor_{(\W\BA)\ot}\ar[dr]_{\Tor_{\e\ot}}
&
&\Tor_{\W\B(A\ot)} \ar[ll]_{\Tor_\psi} \ar[dl]^{\Tor_\e{\vphantom{\Tor_{\e\ot}}}}	\\
&	\Tor_{A\ot}		\mathrlap,	
}
\]
which commutes by \Cref{def:psi} and \Cref{thm:WDg}.
%

\nd 3. The triangle/diamond to its right commutes by \Cref{rmk:qi}.

\nd 4. For the lower-left trapezoid,
select a \DG $A$-module resolution $\defm{P^\bul}$ of $X$
as in \Cref{def:bariso}
and write $\defm{\predbariso}\: P^\bul \lt \B(X,A,A) $
for the map of complexes of \DG $A$-modules over $X$,
so that $\predbariso \ox \id_Y$
induces
$\bariso\:\Tor_A(X,Y)\lt\TOR_A(X,Y)$.
Similarly 
take a proper projective \DG $A\ot$-module resolution $\defm{Q^\bul}$ of $X\ot$
and let $\defm{\predsquism}\: Q^\bul \lt \B(X\ot,A\ot,A\ot)$
be a map of complexes of \DG $A$-modules over $X$,
so that $\predsquism \ox \id_{Y\ot}$
induces $\bariso\:\Tor_{A\ot}(X\ot,Y\ot) \lt \TOR_{A\ot}(X\ot,Y\ot)$.
Recalling the notation  ${*}\: \smash{\BA \xepi \e \kk \xmono \h A}$,
let $\defm\tau$ be the twisting cochain
$t^A \ox {*} + {*} \ox t^A\: \BA \ox \BA \lt A \ox A$.
Then one has a diagram of cochain maps
\[
\resizebox{6.5in}{!}{
\xymatrix@C=4em{
	(P^\bul \ox_A Y)\otimes (P^\bul \ox_A Y)
		\ar[r]_\sim
		\ar[d]^{(f \ox \id) \,\otimes\, (f \ox \id)}
		&
	(\Pb \ox \Pb) \ox_{A\ot} (Y \ox Y)
		\ar[r]
		\ar[d]
	&
	Q^\bul \ox_{A\ot} Y\ot\ar[d]^{\predsquism \ox \id}\\
	\dsp
	(X\ox_{xt^A} \B A \ox_{yt^A} Y)\otimes	(X \ox_{xt^A} \B A \ox_{yt^A} Y)
	\ar[r]^(.5)\sim_(.5){(2\,3\,5\,4)}
		&
		\dsp(X \ox X) \!\! \ox_{\htn x\ot\tau} \! (\BA \ox \BA)
			\!\!\ox_{\htn y\ot\tau} \! (Y \ox Y)
			\ar[r]_{\id\ox\nabla\ox\id}
	&
	\dsp(X \ox X) \!\!\!\!\ox_{x\ot t^{A\ot}} \!\!\!\B(A \ox A) \!\!\!\!\ox_{x\ot t^{A\ot\vphantom{x^{x^x}}}} \!\!\!(Y \ox Y)
	\mathrlap,
}
}
\]
in which the unlabeled
horizontal map exists because $\Pb\ox\Pb$ is projective and $Q^\bul$
exact,
the unlabeled vertical map is defined so as to make the first square commute,
and the second square commutes up to homotopy
by the essential uniqueness of a map from a projective complex to a resolution.
By definition, the composite along the top induces $\ext$ in cohomology,
and the outer vertical maps respectively induce $\bariso \ox \bariso$ and $\bariso$.


\nd 5. The upper-right square commutes by naturality of $\WBbariso$ in \Cref{def:bariso}.

\nd 6. We subdivide a square of bar constructions inducing 
the lower-right square of \eqref{eq:comparison-squares}:
\quation{\label{eq:subdiv}
\begin{aligned}
\mathclap{
\resizebox{6.25in}{!}{
\xymatrix@C=7em@R=3em{
\dsp\W\BX\ot \!\!\! 
	\ox_{t^{\W\BX\ot}\W\B\W x\ot} \!\!\! 
	\B\W\BA\ot \!\!\! \ox_{t^{\W\BY\ot}\W\B\W y\ot} \!\!\! \W\BY\ot 
									\ar[r]^(.55){\W\Phi\ox\W\B\Phi\ox\W\Phi} 
									\ar[d]|(.5)\hole|(.55){\e \ox \id \ox \e}
											|(.6)\hole
&
\dsp\W\BX \!\! \ox_{t^{\W\BX}\W\B\W x\vphantom{\ot}} \!\! 
	\B\W\BA \!\! \ox_{t^{\W\BY}\W\B\W y\vphantom{\ot}} \!\! \W\BY  
									\ar[d]|(.5)\hole|(.55){\e \ox \id \ox \e}
											|(.6)\hole
\\
\dsp \smash{X\ot} \ox_{t^X\B x\ot\B\e} \B\W\smash{\BA\ot} \ox_{t^Y\B y\ot\B\e} \smash{Y\ot}
								\ar[r]^(.55){\B(\Phi,\B\W\Phi,\Phi)} 
									\ar[d]|(.5)\hole|(.55){\id\ox\B\e\ox\id}
											|(.6)\hole
&
{\dsp X \ox_{t^X\B x\B\e\vphantom{\ot}} \B\W\BA \ox_{t^Y\B y\B\e\vphantom{\ot}} Y  }
									\ar[d]|(.5)\hole|(.55){\id\ox\B\e\ox\id}
											|(.6)\hole
\\
\dsp \smash{X\ot} \ox_{t^X\B x\ot} \smash{\BA\ot} \ox_{t^Y\B y\ot} \smash{Y\ot}
								\ar[r]^(.55){\B(\Phi,\Phi,\Phi)} 
&
\dsp X \ox_{t^X\B x} \BA \ox_{t^Y\B y} Y  
\mathrlap;
}
}
}
\end{aligned}
}
but the constituent squares do not actually commute.

We substitute the top square of \eqref{eq:subdiv} with 
\[
\xymatrix@C=7.5em{
	\B\big(\W\B(X\ot),\W\B(A\ot),\W\B(Y\ot)\big) 	
						\ar[r]^(.57){\W\Phi_X \ox \B\W\Phi_A \ox \W\Phi_Y}	&
	 \B(\W\BX,\W\BA,\W\BY) 				\ar[d]^{\e \ox \id \ox \e}		\\
	 \B\big(X\ot,\W\B(A\ot),Y\ot\big) 	\ar[u]^{\B(\h,\,\id,\,\h)}
										\ar[r]_(.57){\B(\Phi_X,\B\W\Phi_A,\Phi_Y)}
																		 &
	 \B(X,\W\BA,Y)
	 \mathrlap,
	}
\]
where $\B(\h,\id,\h)$ is a section of the quasi-isomorphism $\e \ox \id \ox \e$,\footnote{\ 
	To see this, apply $\e$ to the first tensor-factor of $\B(\h,\id,\h)$ 
	as written in the display below,
	recovering $\e t^{\W\BX\ot}\h(\desusp_X \ox \B\e \o \B\W\B x\ot)$.
	By \Cref{thm:twisting-adjunction},
	one has $\e t^{\W\BX\ot}\h = \e t_{\BX\ot} = t^{X\ot}A$,
	which annihilates $\B_{\geq 2} X\ot$.
	Thus for this map to have a nonzero value, the $\B\e \o \B\W\B x\ot$
	must be valued in $\kk = \B_0 X\ot$ rather than $\B_{\geq 1} X\ot$,
	since $\desusp_X$ is already valued in $\B_1 A$.
    This implies the contribution of $\smash{\iter{\D_{\B\W\B(A\ot)}}3}$
	to this tensor factor must be $1$.
	The same holds of the third tensor-factor,
    so the interesting 
    terms of $\smash{\iter{\D_{\B\W\B(A\ot)}}3}$
	are killed and the factor 
    $\big(\!\id \ox \smash{\iter{\D_{\B\W\B(A\ot)}}3} \ox \id\!\big)$
	is functionally just
	$\id_{X\ot} \ox \id_{\B\W\B(A\ot)} \ox \id_{Y\ot}$.
	Looking at the third tensor factor for a change,
	we then have
	$
    \smash{\e\o t^{\W\BY\ot}\o \h\o \desusp_Y} = 
    \smash{\e\o t_{\BY\ot}\o \desusp_Y} = 
    \smash{t^{Y\ot}}\o \desusp_Y = 
    \smash{\id_{\ol{Y\ot}}}
	$,
	again by \Cref{thm:twisting-adjunction}. 
}
and claim this square actually does commute.
Unpacking the definitions~\cite[Prop.~1.26]{carlsonfranzlong},
the composition along the top is
\[
	\begin{multlined}
(\e \ox \id \ox \e) \o (\W\Phi_X \ox \B\W\Phi_A \ox \W\Phi_Y)
\\
\o
\big(
t^{\W\B X\ot}\h(\desusp_X \ox \B\e \o \B\W\B \,\mn x\ot) 
\ox
\id
\ox
t^{\W\B Y\ot}\h(\B\e \o \B\W\B \,\mn y\ot \ox \desusp_Y)
\big) 
\o
\big(\!\id \ox \iter{\D_{\B\W\B(A\ot)}}3 \ox \id\!\big)\mathrlap.
\end{multlined}
\]
and the bottom map 
$\B(\Phi_X,\B\W\Phi_A,\Phi_Y)$ 
is
\eqn{
\big(
	t^X\Phi_X(\desusp_X \ox \B\e \o \B\W\B \,\mn x\ot) 
	&\ox
	\B\W\Phi_A
	\ox
	t^Y\Phi_Y(\B\e \o \B\W\B \,\mn y\ot \ox \desusp_Y)
\big)
	 \o
\big(\!\id \ox \iter{\D_{\B\W\B(A\ot)}}3 \ox \id\!\big)\mathrlap.
}
The initial comultiplication is the same in both cases,
and after, both maps are compositions of tensor products
of maps on the $X$, $A$, and $Y$ components we may examine separately.
It is easy to see that the ``$A$'' map in both cases is $\B\W\Phi_A$.
The $X$ and $Y$ components are symmetrical, and for variety we check
the $Y$ component and suppress the check for $X$ this time.
That the $Y$ components are equal is the claim that
\[
	t^Y\Phi_Y(\B\e \o \B\W\B \,\mn y\ot \ox \desusp_Y)
\ =\ 
\e \o \W\Phi_Y \o t^{\W\B Y\ot}\h(\B\e \o \B\W\B \,\mn y\ot \ox \desusp_Y)
\mathrlap.
\]
It will evidently be enough to check that
\[
t^Y\Phi_Y
\ =\ 
\e \o \W\Phi_Y \o t^{\W\B Y\ot}\h
\mathrlap.
\]
But this, finally, amounts to the commutativity of the following diagram,
in which the parallelogram follows by naturality of $t_{(-)}$
and the triangles by \Cref{thm:twisting-adjunction}.
\[
\begin{aligned}
\xymatrix@C=1.5em@R=3em{
	\B\W\B Y\ot 
				\ar[d]_{t^{\W\B Y\ot}} 					&
	\B Y\ot		\ar[l]_(.4)\h
				\ar[dl]^(.4){t_{\B Y^{\otimes 2}}}
				\ar[r]^(.55){\Phi_Y}							&
	\B Y		\ar[dl]_(.5){t_{\B Y}}
				\ar[d]^{t^Y}							\\
	\W\B Y\ot	\ar[r]_{\W\Phi_Y}						&
	\W\B Y		\ar[r]_(.6){\e}								&
	Y					
	}
\end{aligned}
\]

We substitute the bottom square of \eqref{eq:subdiv} with 
\[
\xymatrix@C=6em{
	 \B\big(X\ot,\W\B(A\ot),Y\ot\big) 	
										\ar[r]^(.575){\B(\Phi_X,\B\W\Phi_A,\Phi_Y)}
																		 &
	 \B(X,\W\BA,Y)						\ar[d]^{\id\ox\B\e\ox\id}		\\
	 \B(X\ot,A\ot,Y\ot)
						\ar[u]^{\id\ox\h\ox\id}
						\ar[r]^(.55){\B(\Phi_X,\Phi_A,\Phi_Y)}	&
	 \B(X,A,Y)
	 \mathrlap,
	}
\]
where again $\id \ox \h \ox \id$ is a section of $\id \ox \B\e \ox \id$
by the unit--counit identities 
    for the $\W \adj \B$ adjunction.
The bottom map $\B(\Phi_X,\Phi_A,\Phi_Y)$ is 
\[
\big(
	t^X\Phi_X(\desusp_X \ox \B \,\mn x\ot) 
	\ox
	\Phi_A
	\ox
	t^Y\Phi_Y(\B \,\mn y\ot \ox \desusp_Y)
\big)
	 \o
\big(\!\id \ox \iter{\D_{\B(A\ot)}}3 \ox \id\!\big)
\]
and the composition along the top is 
\[
\begin{multlined}
(\id\ox\B\e\ox\id)
\big(
	t^X\Phi_X(\desusp_X \ox \B\e \o \B\W\B \,\mn x\ot) 
	\ox
	\B\W\Phi_A
	\ox
	t^Y\Phi_Y(\B\e \o \B\W\B \,\mn y\ot \ox \desusp_Y)
\big)
\\
\big(\!\id \ox \iter{\D_{\B\W\B(A\ot)}}3 \ox \id\!\big)
(\id\ox\h\ox\id)\mathrlap.
\end{multlined}
\]
We have
$\iter{\D_{\B\W\B(A\ot)}}3 \o \h = \h^{\otimes 3} \o\iter{\D_{\B(A\ot)}}3$,
matching the first factor of the bottom map 
and passing an $\h$ forward in each of the three tensor-factors.
The ``$A$'' tensor-factor is
\[
\B\e \o \B\W\Phi_A \o \h = \B\e \o \h \o \Phi_A = \Phi_A
\]
by naturality of $\h$ and the unit--counit identities,
matching that of the bottom map.
For the ``$X$'' factor, one gets 
\[
	t^X\Phi_X(\desusp_X \ox \B\e \o \B\W\B \,\mn x\ot \o \h)
		=
	t^X\Phi_X(\desusp_X \ox \B\e \o \h \o \B \,\mn x\ot)
		=
	t^X\Phi_X(\desusp_X \ox \B \,\mn x\ot)\mathrlap,
\]
again matching the other map, and matching the ``$Y$'' factors is similar.
\end{proof}

\brmk
It thus seems even more plausible than it did in the previous work~\cite[Rmk.~A.26]{carlsonfranzlong}
that the product \eqref{eq:bar-prod} on the two-sided bar construction
is the binary component in a sequence of operations making 
it an \Ai-algebra, but we will not try to prove this here.
\ermk
\bs

{\footnotesize\bibliography{bibshort} }

\begin{thebibliography}{BaMR14}

\bibitem[BaMR14]{barthelmayriehl2014}
Tobias Barthel, J.~Peter May, and Emily Riehl.
\newblock Six model structures for {DG}-modules over {DGAs}: model category
theory in homological action.
\newblock {\em New York J. Math}, 20:1077--1159, 2014.
\newblock \url{http://nyjm.albany.edu/j/2014/20-53.html}, \href
{http://arxiv.org/abs/1310.1159} {\path{arXiv:1310.1159}}, \href
{http://dx.doi.org/10.2140/agt.2013.13.1089}
{\path{doi:10.2140/agt.2013.13.1089}}.

\bibitem[BeF04]{bergerfresse2004operad}
Clemens Berger and Benoit Fresse.
\newblock Combinatorial operad actions on cochains.
\newblock {\em Math. Proc. Cambridge Philos. Soc.}, 137(1):135--174, 2004.
\newblock \href {http://arxiv.org/abs/0109158} {\path{arXiv:0109158}}, \href
{http://dx.doi.org/10.1017/S0305004103007138}
{\path{doi:10.1017/S0305004103007138}}.

\bibitem[Bor51]{borel1951leray}
Armand Borel.
\newblock {\em {Cohomologie des espaces localement compacts d'apr{\'e}s J.
Leray}}, volume~2 of {\em {Lecture Notes in Math.}}
\newblock Springer-Verlag, 1951.
\newblock S{\'e}m. de Top. alg., ETH.

\bibitem[Car22a]{carlson2022munkholm}
Jeffrey~D. Carlson.
\newblock A ring structure on {Tor}.
\newblock 2022.
\newblock  \url{https://jdkcarlson.github.io/Tor.pdf}, \href
  {http://arxiv.org/abs/2306.04860} {\path{arXiv:2306.04860}}.

\bibitem[Car22b]{carlson2022collapses}
Jeffrey~D. Carlson.
\newblock Collapse results for {Eilenberg--Moore} spectral sequences.
\newblock 2022.
\newblock  \url{https://jdkcarlson.github.io/conf.pdf}.

\bibitem[CaF21]{carlsonfranzlong}
Jeffrey D.~Carlson (appendix joint~with Matthias~Franz).
\newblock The cohomology of biquotients via a product on the two-sided bar
construction (expository version).
\newblock 2021.
\newblock \href {http://arxiv.org/abs/2106.02986v1}
{\path{arXiv:2106.02986v1}}.

\bibitem[CarE]{cartaneilenberg}
Henri Cartan and Samuel Eilenberg.
\newblock {\em Homological Algebra}, volume~19 of {\em Princeton Math. Ser.}
\newblock Princeton Univ. Press, 1999 (1956).
\newblock 
\url{http://www.math.stonybrook.edu/~mmovshev/BOOKS/homologicalalgeb033541mbp.pdf}.

\bibitem[Fr20]{franz2019shc}
Matthias Franz.
\newblock Homotopy {Gerstenhaber} algebras are strongly homotopy commutative.
\newblock {\em J. Homotopy Relat. Struct.}, 15(3):557--595, 2020.
\newblock \href {http://arxiv.org/abs/1907.04778} {\path{arXiv:1907.04778}},
  \href {https://doi.org/10.1007/s40062-020-00268-y}
  {\path{doi:10.1007/s40062-020-00268-y}}.

\bibitem[Fr21]{franz2019homogeneous}
Matthias Franz.
\newblock The cohomology rings of homogeneous spaces.
\newblock {\em J. Topol.}, 14(4):1396--1447, 2021.
\newblock \href {http://arxiv.org/abs/1907.04777} {\path{arXiv:1907.04777}},
  \href {https://doi.org/10.1112/topo.12213} {\path{doi:10.1112/topo.12213}}.

\bibitem[GeV95]{gerstenhabervoronov1995}
Murray Gerstenhaber and Alexander~A. Voronov.
\newblock Homotopy {G}-algebras and moduli space operad.
\newblock {\em Int. Math. Res. Not.}, 1995(3):141--153, 02 1995.
\newblock 
\url{https://academic.oup.com/imrn/article-pdf/1995/3/141/6768479/1995-3-141.pdf},
\href {http://arxiv.org/abs/hep-th/9409063} {\path{arXiv:hep-th/9409063}},
\href {http://dx.doi.org/10.1155/S1073792895000110}.

\bibitem[GuM]{gugenheimmay}
Victor~K.A.M. Gugenheim and J.~Peter May.
\newblock {\em On the Theory and Applications of Differential Torsion
Products}, volume 142 of {\em Mem. Amer. Math. Soc.}
\newblock Amer. Math. Soc., 1974.

\bibitem[HMS74]{husemollermoorestasheff1974}
Dale Husemoller, John~C. Moore, and James Stasheff.
\newblock Differential homological algebra and homogeneous spaces.
\newblock {\em J. Pure Appl. Algebra}, 5(2):113--185, 1974.
\newblock URL:
  \url{http://math.mit.edu/~hrm/18.917/husemoller-moore-stasheff.pdf}.

\bibitem[LH]{LH}
Kenji Lef{\`e}vre-Hasegawa.
\newblock {\em Sur les $A_\infty$-cat{\'e}gories}.
\newblock Ph.D. thesis, Universit{\'e} de Paris VII---Denis Diderot, 2002.
\newblock \href {http://arxiv.org/abs/CT/0310337} {\path{arXiv:CT/0310337}}.

\bibitem[Mac]{maclane}
Saunders Mac~Lane.
\newblock {\em Homology}.
\newblock Classics in Mathematics (Grundlehren Math. Wiss. vol.~114). Springer,
2012 (1963, 1975).

\bibitem[McC]{mcclearyspectral}
John McCleary.
\newblock {\em {A User's Guide to Spectral Sequences}}, volume~58 of {\em
  {Cambridge Stud. Adv. Math.}}
\newblock Cambridge Univ. Press, Cambridge, 2001.

\bibitem[McS03]{mccluresmith2003}
James McClure and Jeffrey Smith.
\newblock Multivariable cochain operations and little $n$-cubes.
\newblock {\em J. Amer. Math. Soc.}, 16(3):681--704, 2003.
\newblock \href {http://arxiv.org/abs/0106024} {\path{arXiv:0106024}}, \href
{http://dx.doi.org/10.1090/S0894-0347-03-00419-3}
{\path{doi:10.1090/S0894-0347-03-00419-3}}.

\bibitem[Mun74]{munkholm1974emss}
Hans~J. Munkholm.
\newblock {The {Eilenberg--Moore} spectral sequence and strongly homotopy
  multiplicative maps}.
\newblock {\em J. Pure Appl. Algebra}, 5(1):1--50, 1974.
\newblock \href {https://doi.org/10.1016/0022-4049(74)90002-4}
  {\path{doi:10.1016/0022-4049(74)90002-4}}.

\bibitem[Mun78]{munkholm1978dga}
Hans~J. Munkholm.
\newblock {DGA} algebras as a {Quillen} model category: {Relations} to shm
maps.
\newblock {\em J. Pure Appl. Algebra}, 13(3):221--232, 1978.
\newblock \href {http://dx.doi.org/10.1016/0022-4049(78)90009-9}
{\path{doi:10.1016/0022-4049(78)90009-9}}.

\bibitem[Pr11]{proute2011}
Alain Prout{\'e}.
\newblock A${}_\infty$-structures. {Mod{\`e}les} minimaux de {Baues}--{Lemaire}
et {Kadeishvili} et homologie des fibrations.
\newblock {\em Repr. Theory Appl. Categ}, 21:1--99, 2011.
\newblock  \url{http://tac.mta.ca/tac/reprints/articles/21/tr21abs.html}.

\bibitem[Sm67]{smith1967emss}
Larry Smith.
\newblock Homological algebra and the {Eilenberg}--{Moore} spectral sequence.
\newblock {\em Trans. Amer. Math. Soc.}, 129:58--93, 1967.
\newblock \href {http://dx.doi.org/10.2307/1994364}
{\path{doi:10.2307/1994364}}.

\bibitem[StH70]{halperinstasheff1970}
James Stasheff and Steve Halperin [{\em sic}].
\newblock Differential algebra in its own rite [{\em sic}].
\newblock In {\em Proc. Adv. Study Inst. Alg. Top.(Aarhus 1970)}, volume~3,
pages 567--577, 1970.
  \url{https://u.pcloud.link/publink/show?code=XZW91hVZsnmlRWXLNlVM4sdEW7OE1QsWunx7}.

\bibitem[Wolf77]{wolf1977homogeneous}
Joel~L. Wolf.
\newblock {The cohomology of homogeneous spaces}.
\newblock {\em Amer. J. Math.}, pages 312--340, 1977.
\newblock \href {http://dx.doi.org/10.2307/2373822}
{\path{doi:10.2307/2373822}}.
\end{thebibliography}

\bs

\nd\footnotesize{%
	\url{jeffrey.carlson@tufts.edu
	}
}
\end{document}